\documentclass{amsart}
\usepackage{tikz,hyperref}

\newtheorem{theorem}{Theorem}
\newtheorem{corollary}{Corollary}

\newtheorem{lemma}{Lemma}
\newtheorem{observation}{Observation}
\theoremstyle{definition}
\newtheorem{definition}{Definition}
\newtheorem{example}{Example}
\newtheorem{remark}{Remark}
\numberwithin{equation}{section}


\begin{document}
\title[Critical Itineraries]{Critical Itineraries of Maps with Constant Slope and One Discontinuity}
\author[M. Barnsley]{Michael Barnsley}
\address{The Australian National University \\
Canberra, Australia}
\email{michael.barnsley@anu.edu.au}
\author[W. Steiner]{Wolfgang Steiner}
\address{LIAFA, CNRS UMR 7089, Universit\'e Paris Diderot -- Paris 7 \\ Paris, France}
\email{steiner@liafa.univ-paris-diderot.fr}
\author[A. Vince]{Andrew Vince}
\address{Department of Mathematics, University of Florida \\
Gainesville, FL, USA }
\email{avince@ufl.edu}

\begin{abstract}
For a function from the unit interval to itself with constant slope and one discontinuity, the itineraries of the point of discontinuity are called the critical itineraries. 
These critical itineraries play a significant role in the study of $\beta$-expansions (with positive or negative~$\beta$) and fractal transformations. 
A~combinatorial characterization of the critical itineraries of such functions is provided. 
\end{abstract}
\maketitle

\section{Introduction} \label{sec:Intro}

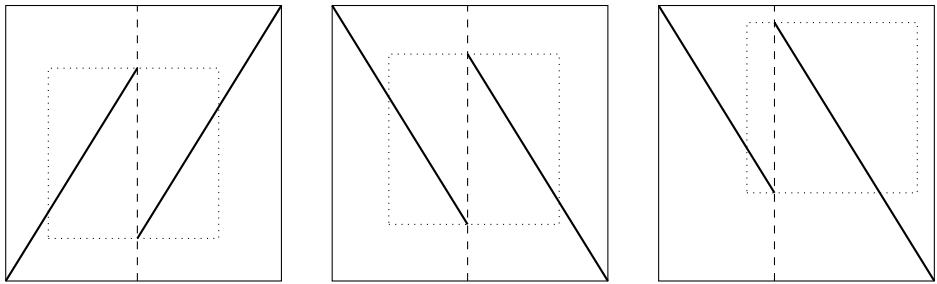
\begin{figure}[ht]
\centerline{\begin{tikzpicture}[scale=1.4]
\draw(0,0)--(2.618,0)--(2.618,2.618)--(0,2.618)--cycle;
\draw[thick](0,0)--(1.25,2.0225) (1.25,.4045)--(2.618,2.618);
\draw[dashed](1.25,0)--(1.25,2.618);
\draw[dotted](.4045,.4045)--(2.0225,.4045)--(2.0225,2.0225)--(.4045,2.0225)--cycle;
\begin{scope}[shift={(4.1,1)}]
\draw(-1,-1)--(1.618,-1)--(1.618,1.618)--(-1,1.618)--cycle;
\draw[thick](-1,1.618)--(.286,-.462) (.286,1.156)--(1.618,-1);
\draw[dashed](.286,-1)--(.286,1.618);
\draw[dotted](-.462,-.462)--(1.156,-.462)--(1.156,1.156)--(-.462,1.156)--cycle;
\end{scope}
\begin{scope}[shift={(7.2,1)}]
\draw(-1,-1)--(1.618,-1)--(1.618,1.618)--(-1,1.618)--cycle;
\draw[thick](-1,1.618)--(.1,-.162) (.1,1.456)--(1.618,-1);
\draw[dashed](.1,-1)--(.1,1.618);
\draw[dotted](-.162,-.162)--(1.456,-.162)--(1.456,1.456)--(-.162,1.456)--cycle;
\end{scope}
\end{tikzpicture}}
\caption{Maps with constant slope and one discontinuity.} \label{fig:Fig1}
\end{figure}

The dynamics of a function from the unit interval to itself is a topic with a long history. 
While most results concern continuous functions, this paper deals with the dynamics of the archetypal families of discontinuous dynamical systems illustrated in Figure~\ref{fig:Fig1}.
These discontinuous functions with constant slope, formally defined below, are often chosen as canonical representatives of conjugacy classes of Lorenz maps \cite{Gl,GS,K}.   
The Lorenz maps serve as models for Poincar\'{e} return maps for Lorenz flows~\cite{L} and play a central role in recent work in fractal geometry~\cite{Ba,BHV}.

Continuous non-differentiable transformations, used in digital imaging and 3D
printing applications, can be constructed using conjugate pairs of such discontinuous systems~\cite{BHI}.   Parameterized families of
such discontinuous systems, and others that are conjugate or semiconjugate to
them, occur in models for a large class of engineering applications
such as circuits, electronics, control systems, and phenomena such as earthquakes; see \cite{HHG,JB,SC}
 and references therein.

\begin{figure}[ht]
\centerline{\begin{tikzpicture}[scale=1.4]
\draw(0,0)--(2.618,0)--(2.618,2.618)--(0,2.618)--cycle;
\draw[thick](0,0)--(1,1.618) (1,0)--(2.618,2.618);
\draw[dashed](1,0)--(1,2.618);
\draw[dotted](1.618,0)--(1.618,1.618)--(0,1.618);
\begin{scope}[shift={(5,1)}]
\draw(-1,-1)--(1.618,-1)--(1.618,1.618)--(-1,1.618)--cycle;
\draw[thick](-1,1.618)--(.382,-.618) (.382,1)--(1.618,-1);
\draw[dashed](.382,-1)--(.382,1.618);
\draw[dotted](-.618,-.618)--(1,-.618)--(1,1)--(-.618,1)--cycle;
\end{scope}
\end{tikzpicture}}
\caption{Greedy $\beta$-transformation and Ito-Sadahiro's ($-\beta$)-transformation.} \label{fig:Fig2}
\end{figure}
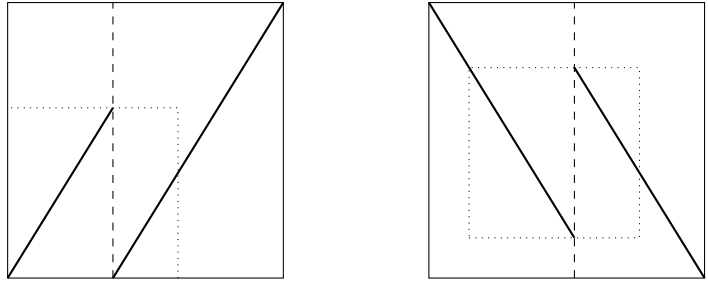

The dynamics of $\beta$-transformations --- functions of the type depicted in Figure~\ref{fig:Fig2} restricted to the inner square --- are integral to the study of the representation of the real numbers using non-integer bases.  
For positive~$\beta$, there is a large literature on this subject beginning with the pioneering work of R\'enyi and Parry \cite{R,P1}.
Generalizations such as linear mod one functions --- depicted in the left and middle panels of Figure~\ref{fig:Fig1} --- have also been studied extensively; see for example \cite{H2,FL}.
The study of negative $\beta$-transformations, which were often neglected, gained a new momentum with the paper~\cite{IS} by Ito and Sadahiro. 
Many arguments for positive slopes easily adapt to negative ones, but some properties of the positive case are not true for the negative slopes, see e.g.\ \cite{LS} and Example~\ref{ex:3} in Section~\ref{sec:main-results}.

 It is well known that the behavior of such discontinuous dynamical systems is mediated by the critical itineraries, namely certain symbolic orbits that are defined below. 
A~similar situation occurs for the dynamics of continuous systems, for which the canonical representative is the family of logistic maps $L_{a}: [0,1] \rightarrow [0,1]$, where $L_{a}(x)=ax(1-x)$, $a\in(0,4]$; see for example \cite{CE}. In
this case, conditions under which a given continuous system is conjugate to a
logistic map are well understood in terms of a symbolic orbit of the critical
point, $x=0.5$.  This symbolic orbit (which is analogous to but not the same
as the critical itineraries in the present work) have been fully characterized~\cite{MT}.
The present paper provides an analogous, succinct, complete characterization of
the critical itineraries of the discontinuous systems illustrated in Figure 1.
While related results are present in the literature, many discussed in Section~\ref{sec:G},  the present characterization appears new; in particular, the case of negative slopes has not been treated elsewhere.

In the left and middle images in Figure~\ref{fig:Fig1}, the restriction to the dotted square is (after proper renormalization) of the form $\beta x + \alpha \bmod 1$ with $|\beta| > 1$. 
This is true in general when $\beta > 1$. 
For $\beta < -1$, however, the right image in Figure~\ref{fig:Fig1} gives an example where the situation is different, hence the class of functions that we consider is larger than that of the maps $\beta x + \alpha \bmod 1$ with one discontinuity.

It is convenient to consider \textbf{generalized $\beta$-transformations} of the form
\[
f_{\beta,p}:\ \mathbb{R} \to \mathbb{R}, \quad x \mapsto \begin{cases} \beta x & \text{if}\ x < p, \\[.5ex] \beta x\ \text{or}\ \beta (x-1) & \text{if}\ x = p, \\[.5ex] \beta (x - 1) & \text{if}\ x > p,\end{cases}
\]
with $\beta, p \in \mathbb{R}$, and $|\beta| > 1$.
More precisely, we define two functions $f_{\beta,p,\pm}$ by 
\[
f_{\beta,p,-}(p) = \beta p,\ f_{\beta,p,+}(p) = \beta (p-1),\ f_{\beta,p,-}(x) = f_{\beta,p,+}(x) = f_{\beta,p}(x)\ \text{for}\ x \ne p.
\]
For the trajectories $f_{\beta,p,\pm}^n(p)$ of the discontinuity to be bounded, we need that
\begin{equation} \label{e:betap}
\beta > 1,\ 1 \le p \le \frac{1}{\beta-1}, \quad \text{or} \quad \beta < -1,\ \frac{\beta^2+\beta-1}{\beta^2-1} \le p \le \frac{1}{\beta^2-1}
\end{equation}
(which implies that $|\beta| \le 2$). 
For these parameters, we have
\[
f_{\beta,p}\big(\big[0, \tfrac{\beta}{\beta-1}\big]\big) = \big[0, \tfrac{\beta}{\beta-1}\big] \quad \text{and} \quad f_{\beta,p}\big(\big[\tfrac{\beta}{\beta^2-1}, \tfrac{\beta^2}{\beta^2-1}\big]\big) = \big[\tfrac{\beta}{\beta^2-1}, \tfrac{\beta^2}{\beta^2-1}\big]
\]
when $\beta > 1$ and $\beta < -1$, respectively.
The restriction of~$f_{\beta,p}$ to the respective interval has the form of a map in Figure~\ref{fig:Fig1}.
Moreover, every expanding map from the unit interval to itself with constant slope and one discontinuity is conjugate to the restriction to some interval of some function $f_{\beta,p,-}$ or~$f_{\beta,p,+}$.

The trajectories of points in~$\mathbb{R}$ by $f_{\beta,p,\pm}$ can be coded by elements of 
\[
\Omega= \{0,1\}^\omega,
\] 
which denotes the set of infinite words (or sequences) $\mathbf{c} = c_0 c_{1} c_{2} \cdots$ on the alphabet $\{0,1\}$. 
For $x \in \mathbb{R}$, the two \textbf{itineraries} of~$x$ are 
\begin{align*}
\tau_{\beta,p,-}(x) & = c_0 c_1 \cdots \quad \text{with} \quad c_n = \begin{cases}0 & \text{if} \ f_{\beta,p,-}^n(x) \le p, \\[.5ex] 1 & \text{if}\ f_{\beta,p,-}^n(x) > p,\end{cases} \\
\tau_{\beta,p,+}(x) & = c_0 c_1 \cdots \quad \text{with} \quad c_n = \begin{cases}0 & \text{if} \ f_{\beta,p,+}^n(x) < p, \\[.5ex] 1 & \text{if}\ f_{\beta,p,+}^n(x) \ge p.\end{cases} 
\end{align*}
The two itineraries of the point of discontinuity~$p$ play a special role. 
Call $\tau_- := \tau_{\beta,p,-}(p)$ and $\tau_+ := \tau_{\beta,p,+}(p)$ the \textbf{critical itineraries} of~$f_{\beta,p}$. 
The pair $(\tau_-; \tau_+)$ is also referred to as the \textit{kneading invariant}
of~$f_{\beta,p}$. 
For $\beta > 1$, the critical itineraries are equal to the \textit{limit itineraries} $\lim_{x\uparrow p} \tau_{\beta,p,\pm}(x)$ and $\lim_{x\downarrow p} \tau_{\beta,p,\pm}(x)$.
For $\beta < -1$, this relation is not necessarily true, see Observation~\ref{o:critneg} in Section~\ref{sec:main-results}.

The main result in this paper is a combinatorial characterization of the critical itineraries of a function~$f_{\beta,p}$.
The possible pairs $(\tau_-; \tau_+)$ are exactly those which are lex-admissible or alt-admissible, as defined in Section~\ref{sec:admissible}.
As a corollary to the main result, we get a characterization of the critical itineraries of $\beta x + \alpha \bmod 1$ (when this map has only one discontinuity). 

For the particular case $p = 1$, $1 < \beta \le 2$, the critical itineraries were already described in~\cite{P1}.
Indeed, we have $\beta\, T_\beta(x) = f_{\beta,1,+}(\beta x)$ for all $x \in [0,1)$, where $T_\beta$ is the greedy $\beta$-transformation, defined by $T_\beta(x) := \beta x - \lfloor \beta x \rfloor$; see also Figure~\ref{fig:Fig2}.
Here, since $\tau_{\beta,1,+}(1) = 1000\cdots$, it is sufficient to study $\tau_{\beta,1,-}(1) = 0\, \tau_{\beta,1,-}(\beta)$.

For $p = \frac{1}{1-\beta}$, $-2 < \beta < -1$, we have $\beta\, T_\beta(x) = f_{\beta,p,+}(\beta x)$ for all $x \in \big[\frac{\beta}{1-\beta}, \frac{1}{1-\beta}\big)$, where $T_\beta$ is the $\beta$-transformation defined in~\cite{IS} by $T_\beta(x) := \beta x - \lfloor \beta x - \frac{\beta}{1-\beta} \rfloor$ (for negative~$\beta$).
The critical itineraries of these maps were characterized in~\cite{St}.
Here, we have $\tau_{\beta,p,-}(p) = 00\, \tau_{\beta,p,-}(\beta^2 p)$ and $\tau_{\beta,p,+}(p) = 1\, \tau_{\beta,p,+}(\beta^2 p)$.

\section{Admissible pairs of words} \label{sec:admissible}
The \textbf{lexicographic} order on~$\Omega$ is the total order defined by $\mathbf{a} < \mathbf{b}$ if $\mathbf{a} \neq \mathbf{b}$ and $a_n < b_n$ where $n$ is the least index such that $a_n \neq b_n$. 
The \textbf{alternating lexicographic} order on~$\Omega$ is the total order defined by $\mathbf{a} < \mathbf{b}$ if $\mathbf{a} \neq \mathbf{b}$ and $(-1)^n\, (a_n-b_n) < 0$, where $n$ is the least index such that $a_n \neq b_n$ (with $\mathbf{a} = a_0 a_1 \cdots$, $\mathbf{b} = b_0 b_1 \cdots$).
We use the notation
\[
(\mathbf{a}, \mathbf{b}) := \{\mathbf{c} \in \Omega :\, \mathbf{a} < \mathbf{c} < \mathbf{b}\}
\]
for the open interval with respect to the specified order (lexicographic or alternating lexicographic); likewise for the closed and half open intervals.

Let $S$ denote the shift operator on $\Omega$, i.e., $S(c_0 c_1 c_2 \cdots) = c_1 c_2 c_3 \cdots$. 
For a set $X \subseteq \Omega$, let
\[
\Omega^X := \big\{\mathbf{c} \in \Omega:\, S^n(\mathbf{c}) \not\in X\ \text{for all}\ n \ge 0\big\}.
\] 
Note that $\Omega^X $ is \emph{shift invariant} in the sense that $S(\Omega^X ) = \Omega^X$.

For $\Lambda\subseteq\Omega$, let 
\[
\Lambda_n := \big\{\mathbf{u} \in \{0,1\}^n:\, \mathbf{u} \mathbf{c} \in \Lambda \;  \text{for some}\ \mathbf{c} \in \Omega \big\}
\]
be the set of length $n$ prefixes of words in~$\Lambda$, and let $|\Lambda_n|$ denote the cardinality of~$\Lambda_n$.  
The \textbf{exponential growth rate} $g(\Lambda)$ of $\Lambda
\subseteq\Omega$ is given by
\[
g(\Lambda):=\lim_{n\rightarrow\infty}\sqrt[n]{|\Lambda_n|},
\]
if the limit exits.  In particular, the limit exists for
$\Lambda =  \Omega^{(\mathbf{a}, \mathbf{b})}$ (see e.g.~\cite{LM}), and 
\[
h\big(\Omega^{(\mathbf{a},\mathbf{b}}\big) :=  \log g\big(\Omega^{(\mathbf{a},\mathbf{b})}\big)
\]
is the topological entropy of $\Omega^{(\mathbf{a},\mathbf{b})}$ considered as a symbolic dynamical system with the shift map~$S$ operating on it (and $\Omega$ equipped with the product topology of the discrete toplogy).  

Let $A := \{0,1\}$. 
The set of finite words over~$A$ is denoted by $A^* = \bigcup_{n\ge0} A^n$. 
The length of a word $\mathbf{u} \in A^*$ is denoted by~$|\mathbf{u}|$, i.e., $|\mathbf{u}| = n$ if $\mathbf{u} \in A^n$. 
The Kleene star $B^* = \bigcup_{n\ge0} B^n$ is also used for sets of words $B \subset A^*$; it denotes the set of finite concatenations of elements from~$B$.
The set of infinite concatenations of elements from~$B$ is denoted by~$B^\omega$. 
For $\mathbf{u} \in A^*$ with $|\mathbf{u}| \ge 1$, the only element of $\{\mathbf{u}\}^\omega$ is the periodic infinite word $\overline{\mathbf{u}} := \mathbf{u} \mathbf{u} \mathbf{u} \cdots \in \Omega$.

\begin{definition} \label{def:admissible} 
Call a pair of words $(\mathbf{a}; \mathbf{b})$ with $\mathbf{a} \in 0\, \Omega$, $\mathbf{b} \in 1\, \Omega$, \textbf{lex-admissible} if the properties (\ref{i:adm1})--(\ref{i:adm3}) below hold for the lexicographic order, \textbf{alt-admissible} if (\ref{i:adm1})--(\ref{i:adm3}) hold for the alternating lexicographic order.
\begin{enumerate}
\itemsep1ex
\item \label{i:adm1}
$S^n(\mathbf{a}) \notin (\mathbf{a}, \mathbf{b}]$ and $S^n(\mathbf{b}) \notin [\mathbf{a}, \mathbf{b})$ for all $n \ge 0$, i.e., $\mathbf{a} \in \Omega^{(\mathbf{a},\mathbf{b}]}$, $\mathbf{b} \in \Omega^{[\mathbf{a},\mathbf{b})}$,
\item \label{i:adm2}
$g\big (\Omega^{(\mathbf{a},\mathbf{b})}\big ) > 1$,
\item \label{i:adm3}
if $\mathbf{a}, \mathbf{b} \in \{\mathbf{u}, \mathbf{v}\}^\omega$ for some finite words $\mathbf{u} \in 0\, \{0,1\}^*$, $\mathbf{v} \in 1\, \{0,1\}^*$, with $\overline{\mathbf{u}} \in \Omega^{(\overline{\mathbf{u}},\overline{\mathbf{v}}]}$, $\overline{\mathbf{v}} \in \Omega^{[\overline{\mathbf{u}},\overline{\mathbf{v}})}$, and $g\big(\Omega^{(\overline{\mathbf{u}},\overline{\mathbf{v}})}\big) = g\big(\Omega^{(\mathbf{a},\mathbf{b})}\big)$, then $\mathbf{a} = \overline{\mathbf{u}}$ and $\mathbf{b} = \overline{\mathbf{v}}$.
\end{enumerate}
\end{definition}

\begin{example}[Pairs with zero exponential growth rate]
It is not hard to find examples of $\mathbf{a}, \mathbf{b} \in \Omega$ satisfying condition~(\ref{i:adm1}) but $g\big(\Omega^{(\mathbf{a},\mathbf{b})}\big) = 1$. 
There are trivial examples such as $\mathbf{a} = \overline{0}$ with arbitrary~$\mathbf{b}$ satisfying~(\ref{i:adm1}), and this is also the case when $\mathbf{a}, \mathbf{b}$ are the critical itineraries of the function $x + \alpha \bmod 1$ with irrational~$\alpha$.  
\end{example}

We define the value of a sequence $\mathbf{c} = c_0 c_1 c_2 \cdots \in \Omega$ in base~$\beta$ by
\[
\langle \mathbf{c} \rangle_\beta := \sum_{n=0}^\infty \frac{c_n}{\beta^n}.
\] 

\section{Main results} \label{sec:main-results}
Our main result is the following theorem, which is proved in Section~\ref{sec:admissible-pair-=}.

\begin{theorem} \label{t:main}
Two words $\mathbf{a}, \mathbf{b} \in \Omega$ are the critical itineraries of~$f_{\beta,p}$ for some $\beta > 1$, $1 \le p \le \frac{1}{\beta-1}$, if and only if the pair $(\mathbf{a}; \mathbf{b})$ is lex-admissible. 

Two words $\mathbf{a}, \mathbf{b} \in \Omega$ are the critical itineraries of~$f_{\beta,p}$ for some $\beta < -1$, $\frac{\beta^2+\beta-1}{\beta^2-1} \le p \le \frac{1}{\beta^2-1}$, if and only if the pair $(\mathbf{a}; \mathbf{b})$ is alt-admissible. 

In either case, we have $|\beta| = g\big(\Omega^{(\mathbf{a},\mathbf{b})}\big)$, $p = \langle \mathbf{a} \rangle_\beta = \langle \mathbf{b} \rangle_\beta$, and $\langle \mathbf{a} \rangle_\gamma \ne \langle \mathbf{b} \rangle_\gamma$ for all $\gamma \in \mathbb{R}$ with $\mathrm{sgn}(\gamma) = \mathrm{sgn}(\beta)$ and $|\gamma| > |\beta|$.
\end{theorem}

The following theorem, which is proved in Section~\ref{sec:expon-growth-rates}, shows that the conditions for $\beta > 1$ can be simplified. 
In fact, the equality $g\big(\Omega^{(\overline{\mathbf{u}},\overline{\mathbf{v}})}\big) = g\big(\Omega^{(\mathbf{a},\mathbf{b})}\big)$ in~(\ref{i:adm3}) is automatically satisfied for the lexicographic order (except for $\overline{\mathbf{u}} = \overline{0}$ or $\overline{\mathbf{v}}  = \overline{1}$, where $g\big(\Omega^{(\overline{\mathbf{u}},\overline{\mathbf{v}})}\big) = 1$). 

\begin{theorem} \label{t:pos}
A~pair of words $(\mathbf{a}; \mathbf{b})$ with $\mathbf{a} \in 0\, \Omega$, $\mathbf{b} \in 1\, \Omega$, is lex-admissible if and only if properties (\ref{i:adm1}) and~(\ref{i:adm2}) of Definition~\ref{def:admissible} and property (\ref{i:pos3}) below hold for the lexicographic order.
\begin{enumerate}
\renewcommand{\theenumi}{3'}
\item \label{i:pos3}
If $\mathbf{a}, \mathbf{b} \in \{\mathbf{u}, \mathbf{v}\}^\omega$ for some finite words $\mathbf{u} \in 01\, \{0,1\}^*$, $\mathbf{v} \in 10\, \{0,1\}^*$, with $\overline{\mathbf{u}} \in \Omega^{(\overline{\mathbf{u}},\overline{\mathbf{v}}]}$ and  $\overline{\mathbf{v}} \in \Omega^{[\overline{\mathbf{u}},\overline{\mathbf{v}})}$, then $\mathbf{a} = \overline{\mathbf{u}}$ and $\mathbf{b} = \overline{\mathbf{v}}$.
\end{enumerate}
\end{theorem}

The following corollary of Theorem~\ref{t:main} characterizes the critical itineraries of ${\beta x + \alpha} \bmod 1$, i.e., the itineraries of the discontinuity point of the maps
\begin{align*}
T_{\beta,\alpha,+}:\ & [0,1) \to [0,1), \quad x \mapsto \beta x + \alpha - \lfloor \beta x + \alpha \rfloor,  \\
T_{\beta,\alpha,-}:\ & (0,1] \to (0,1], \quad x \mapsto \beta x + \alpha - \lceil \beta x + \alpha \rceil + 1, 
\end{align*}
with $\beta > 1$, $0 \le \alpha \le  2-\beta$, or $\beta < -1$, $-\beta - 1 < \alpha < 1$.
For these parameters, both maps $T_{\beta,\alpha,-}$ and $T_{\beta,\alpha,+}$ have a unique discontinuity point, which is at $({1\!-\!\alpha})/\beta$ when $\beta > 1$ and at $-\alpha/\beta$ when $\beta < -1$.
We define the itinerary of $x \in [0,1)$ under~$T_{\beta,\alpha,+}$ as $c_0 c_1 \cdots \in \Omega$ with $c_n = \mathrm{sgn}(\beta)\, \lfloor \beta\, T_{\beta,\alpha,+}^n(x) + \alpha \rfloor$, and the itinerary of $x \in (0,1]$ under~$T_{\beta,\alpha,-}$ as $c_0 c_1 \cdots \in \Omega$ with $c_n = \mathrm{sgn}(\beta)\, (\lceil \beta\, T_{\beta,\alpha,-}^n(x) + \alpha \rceil - 1)$.

\begin{corollary} \label{c:main}
Two words $\mathbf{a}, \mathbf{b} \in \Omega$ are the critical itineraries of $\beta x + \alpha \bmod 1$ for some $\beta > 1$, $0 \le \alpha \le 2-\beta$, if and only if the pair $(\mathbf{a}; \mathbf{b})$ is lex-admissible. 
In this case, we have $\beta = g\big(\Omega^{(\mathbf{a},\mathbf{b})}\big)$ and $\alpha = (\langle \mathbf{a} \rangle_\beta - 1)\, (\beta-1) = (\langle \mathbf{b} \rangle_\beta - 1)\, (\beta-1)$.

Two words $\mathbf{a}, \mathbf{b} \in \Omega$ are the critical itineraries of $\beta x + \alpha \bmod 1$ for some $\beta < -1$, $-\beta - 1 < \alpha < 1$, if and only if the pair $(\mathbf{a}; \mathbf{b})$ is alt-admissible, $S^2(\mathbf{a}) < S(\mathbf{b})$, and $S^2(\mathbf{b}) > S(\mathbf{a})$ (in the alternating lexicographic order). 
In this case, we have $\beta = -g\big(\Omega^{(\mathbf{a},\mathbf{b})}\big)$ and $\alpha = \langle \mathbf{a} \rangle_\beta\, (1-\beta) = \langle \mathbf{b} \rangle_\beta\, (1-\beta)$.
\end{corollary}

As mentioned in the Introduction, what we call critical itineraries are not necessarily the limit itineraries from the left and right to~$p$, when $\beta$ is negative. 
The relation between these two notions is described by the following observation.

\begin{observation} \label{o:critneg}
Let $\beta < -1$, $\frac{\beta^2+\beta-1}{\beta^2-1} \le p \le \frac{1}{\beta^2-1}$. If $\mathbf{a} := \tau_{\beta,p,-}(p)$ is periodic with odd period length, let $\mathbf{u}$ be its primitive period.  
If $\mathbf{b} := \tau_{\beta,p,+}(p)$ is periodic with odd period length, let $\mathbf{v}$ be its primitive period.  
Then
\begin{align*}
\lim_{x\uparrow p} \tau_{\beta,p}(x) & = \begin{cases} 
\mathbf{a} & \text{if $\mathbf{a}$ is not periodic with odd period length}, \\ 
\mathbf{ub} & \text{if $\mathbf{a}$ is periodic with odd period length, but $\mathbf{b}$ is not}, \\ 
\overline{\mathbf{uv}} & \text{if $\mathbf{a}$ and $\mathbf{b}$ are periodic with odd period length}, \end{cases} \\
\lim_{x\downarrow p} \tau_{\beta,p}(x) & = \begin{cases} 
\mathbf{b} & \text{if $\mathbf{b}$ is not periodic with odd period length}, \\ \mathbf{va} & \text{if $\mathbf{b}$ is periodic with odd period length, but $\mathbf{a}$ is not}, \\ 
\overline{\mathbf{vu}} & \text{if $\mathbf{a}$ and $\mathbf{b}$ are periodic with odd period length}.
\end{cases} 
\end{align*}
\end{observation}

\begin{remark}
For the greedy $\beta$-transformation, we have the following:
The pairs of critical itineraries of~$f_{\beta,1}$, $1 < \beta \le 2$, are exactly the pairs $(0 \mathbf{c}; 1 \overline{0})$ with $\mathbf{c} \in 1\, \Omega$ and $1 \overline{0} \ne S^n(\mathbf{c}) \le S(\mathbf{c})$ for all $n \ge 1$ (w.r.t.\ the lexicographic order); cf.~\cite{P1}.

For the cases corresponding to Ito-Sadahiro's $(-\beta)$-transformations, we obtain the following characterization from~\cite{St}.
Let $\mathbf{d} = 100111001001001110011 \cdots \in \Omega$ be the word starting with $\varphi^n(1)$ for all $n \ge 0$, where $\varphi$ denotes the morphism on $\{0,1\}^*$ defined by $\varphi(1) = 100$, $\varphi(0) = 1$. 
Then the critical itineraries of~$f_{\beta,1/(1-\beta)}$, $-2 \le \beta < -1$, are exactly the pairs $(00\mathbf{c}; 1\mathbf{c})$ with $\mathbf{c} \in 1\, \Omega$ such that $S^n(\mathbf{c}) \le \mathbf{c}$ for all $n \ge 1$,  $\mathbf{c} > \mathbf{d}$, and $\mathbf{c} \notin \{\mathbf{u}00, \mathbf{u}1\}^\omega$ for all $\mathbf{u} \in \{0,1\}^*$ with  $\overline{\mathbf{u}1} > \mathbf{d}$, and the pairs $(\overline{00\mathbf{w}}; \overline{1\mathbf{w}})$ with $\mathbf{w} \in 1\, \{0,1\}^*$ such that $S^n(\overline{\mathbf{w}1}) \le \overline{\mathbf{w}1}$ for all $n \ge 1$, $\overline{\mathbf{w}1} > \mathbf{d}$, and $\mathbf{w}1 \notin \{\mathbf{u}00, \mathbf{u}1\}^*$ for all $\mathbf{u} \in \{0,1\}^* \setminus \{\mathbf{w}\}$ with  $\overline{\mathbf{u}1} > \mathbf{d}$, where the inequalities refer to the alternating lexicographic order.

It should be mentioned that $f_{-2,1/3,+}$ is not conjugate to the map~$T_{-2}$ from \cite{IS,St} because $f_{-2,1/3,+}(4/3) = -2/3 \ne 4/3 = -2\, T_{ -2}(-2/3)$. 
The critical itineraries of $f_{-2,1/3}$ are $(00\overline{10}; 00\overline{10})$, which is a pair satisfying the conditions above. 
However, the word~$\overline{10}$ does not satisfy condition~(1.8) in~\cite{St} because $\overline{10} \in \{2, 10\}^\omega$. 

Note also that the inequality $\mathbf{c} > \mathbf{d}$ implies that $g\big(\Omega^{(00\mathbf{c},1\mathbf{c})}\big) > 1$.
Moreover, it is not necessary to verify the equation $g\big(\Omega^{(00\mathbf{c},1\mathbf{c})}\big) = g\big(\Omega^{(\overline{00\mathbf{u}},\overline{1\mathbf{u}})}\big)$ here.
\end{remark}

\begin{example}[Primality Tester] \label{ex:prime} 
The pair 
\[
\begin{aligned} \mathbf a & = 0\,1\,1\,0\,1\,0 \,1\,0\,0\,0\,1\,0\,1 \cdots \\ \mathbf b & = 1\,0\, 0\,0 \cdots, \end{aligned}
\]
where $a_{n} = 1$ if and only if $n+1$ is prime is lex-admissible.
It is easy to check conditions (\ref{i:adm1}) and~(\ref{i:adm3}) of Definition~\ref{def:admissible}. 
Concerning condition (\ref{i:adm2}), we have $\Omega^{(\mathbf{a},\mathbf{b})} \supset \Omega^{(\overline{01},1\overline{0})}$ and thus $g\big(\Omega^{(\mathbf{a},\mathbf{b})}\big) \ge g\big(\Omega^{(\overline{01},1\overline{0})}\big) = (\sqrt{5}+1)/2$.
Therefore, $\mathbf{a}$ and $\mathbf{b}$ are the critical itineraries of $f_{\beta,1}$, where $\beta \approx 1.79$, as stated in Theorem~\ref{t:main}.
By the definition of the critical itineraries, the natural number $n$ is prime if and only if $f_{\beta,1,+}^{n-1}(p) > p$.
In other words, to test whether $n+1$ is prime, we apply the $n^{th}$ iterate of $f_{\beta,1,+}$ to the point of discontinuity~$1$. 
If this iterate lies to the right of~$p$, then $n$ is prime; otherwise, it is composite. 
Two comments are in order.
First, this result has little to do with number theory. 
Second, the method is numerically problematic because~$\beta$, being an irrational number, can be estimated to at most finitely many places.
\end{example}

\begin{example}[Words $\mathbf{a}, \mathbf{b} \in \{\mathbf{u}, \mathbf{v}\}^\omega$ with $g\big(\Omega^{(\mathbf{a},\mathbf{b})}\big) > g\big(\Omega^{(\overline{\mathbf{u}},\overline{\mathbf{v}})}\big) > 1$ for the alternating lexicographic order] \label{ex:3}
The following example illustrates that $\langle \mathbf{a} \rangle_\beta = \langle \mathbf{b} \rangle_\beta$ does not necessarily have a unique solution and that it might be difficult to avoid the condition $g(\Omega^{(\overline{\mathbf{u}},\overline{\mathbf{v}})}) = g(\Omega^{(\mathbf{a},\mathbf{b})})$ in property~(\ref{i:adm3}) of Definition~\ref{def:admissible} for negative~$\beta$. 
Let
\[
\mathbf{u} = 001100000, \quad \mathbf{v} = 110, \quad \mathbf{a} = \mathbf{u} \overline{\mathbf{uv}}, \quad \mathbf{b} = \mathbf{v} \overline{\mathbf{vu}}.
\]
By Lemma~\ref{l:K} below, we can calculate
\[
g\big(\Omega^{(\mathbf{a},\mathbf{b})}\big) = -\beta > -\gamma = g\big(\Omega^{(\overline{\mathbf{u}},\overline{\mathbf{v}})}\big),
\]
with $\beta \approx -1.135888346$ satisfying $\beta^9 = - \beta^6 - 1$, and $\gamma \approx -1.123732821$ satisfying $\gamma^5 =\gamma^4 - \gamma^2 + \gamma -1$.
Note that $1/\gamma$ is a root of the power series $K(z)$ defined by $\mathbf{a}$ and~$\mathbf{b}$ in Lemma~\ref{l:K}, but the largest negative root is~$1/\beta$. 
Setting
\begin{align*}
p & := \langle \mathbf{a} \rangle_\beta = \langle \mathbf{b} \rangle_\beta = \frac{\beta^8+\beta^5+\beta^4-\beta^2-1}{\beta^{12}-1} \approx 0.070528093, \\
q & := \langle \overline{\mathbf{u}} \rangle_\gamma = \langle \overline{\mathbf{v}} \rangle_\gamma = \langle \mathbf{a} \rangle_\gamma = \langle \mathbf{b} \rangle_\gamma = \frac{\gamma^3+\gamma^2}{\gamma^3-1} \approx 0.064590878,
\end{align*}
we have 
\[
\tau_{\beta,p,-}(p) = \mathbf{a}, \quad \tau_{\beta,p,+}(p) = \mathbf{b}, \quad \tau_{\gamma,q,-}(q) = \overline{\mathbf{u}}, \quad \tau_{\gamma,q,+}(q) = \overline{\mathbf{v}}.
\]
According to Corollary~\ref{c:main}, $\mathbf{a}$ and~$\mathbf{b}$ are the critical itineraries of $\beta x + p (1-\beta) \bmod 1$, while $\overline{\mathbf{u}}$ and~$\overline{\mathbf{v}}$ are the critical itineraries of $\gamma x + q (1-\gamma) \bmod 1$.

In this example, we have $g\big(\Omega^{(\mathbf{a},\mathbf{b})}\big) = g\big(\bigcup_{n=0}^\infty S^n \{\mathbf{u}, \mathbf{v}\}^\omega) > g\big(\Omega^{(\overline{\mathbf{u}},\overline{\mathbf{v}})}\big) > 1$.
Such a situation cannot occur for positive~$\beta$, where $g\big(\bigcup_{n=0}^\infty S^n \{\mathbf{u}, \mathbf{v}\}^\omega\big) \le g\big(\Omega^{(\overline{\mathbf{u}},\overline{\mathbf{v}})}\big)$ always holds for $\mathbf{u}, \mathbf{v}$ as in (\ref{i:pos3}); see the proof of Lemma~\ref{l:entpos}.
\end{example}

\section{Admissible pair $=$ critical itineraries}
\label{sec:admissible-pair-=}
\subsection{Address space}
The description of the address space of $f_{\beta,p}$ is fairly standard for positive~$\beta$, see e.g.\ \cite[Theorem~2.5]{KS} or \cite[Theorem~5.1]{BHV}.
We include its proof for completeness and to prepare the slightly more complicated case of negative~$\beta$. 

\begin{lemma} \label{l:addpos}
Let $\beta > 1$, $1 \le p \le \frac{1}{\beta-1}$, $\mathbf{a} := \tau_{\beta,p,-}(p)$, and $\mathbf{b} := \tau_{\beta,p,+}(p)$.
Then the address spaces of $f_{\beta,p,\pm}$ are
\[
\tau_{\beta,p,-}(\mathbb{R}) = \tau_{\beta,p,-}\big(\big[0, \tfrac{\beta}{\beta-1}\big]\big) = 
\Omega^{(\mathbf{a},\mathbf{b}]} \ \text{and} \ \tau_{\beta,p,+}(\mathbb{R}) = 
\tau_{\beta,p,+}\big(\big[0, \tfrac{\beta}{\beta-1}\big]\big) = \Omega^{[\mathbf{a},\mathbf{b})},
\]
with the lexicographic order on~$\Omega$.
In particular, we have $\mathbf{a} \in \Omega^{(\mathbf{a},\mathbf{b}]}$ and $\mathbf{b} \in \Omega^{[\mathbf{a},\mathbf{b})}$.
Moreover, we have $g\big(\Omega^{(\mathbf{a},\mathbf{b})}\big) = \beta$. 
\end{lemma}

\begin{proof}
Let $x \in \big[0, \frac{\beta}{\beta-1}]$, and recall that $f_{\beta,p}\big(\big[0, \tfrac{\beta}{\beta-1}\big]\big) = \big[0, \frac{\beta}{\beta-1}\big]$ for $\beta > 1$, $1 \le p \le \frac{1}{\beta-1}$.
If $\tau_{\beta,p,-}(x)$ agrees with $\mathbf{a}$ or~$\mathbf{b}$ on the first $n$ letters, then $f_{\beta,p,-}^n(x) - f_{\beta,p,\pm}^n(p) = \beta^n (x-p)$.
Therefore, since $f_{\beta,p}^n(x)$ is bounded, $x < p$ implies that $f_{\beta,p,-}^n(x) \le p < f_{\beta,p,-}^n(p)$ for some $n \ge 1$, and $p < x$ implies that $f_{\beta,p,+}^n(p) < p < f_{\beta,p,-}^n(x)$ for some $n \ge 1$, hence $\tau_{\beta,p,-}(x) \le \mathbf{a}$ if $x \le p$ and $\tau_{\beta,p,-}(x) > \mathbf{b}$ if $x > p$.
This gives that $S^n(\tau_{\beta,p,-}(x)) = \tau_{\beta,p,-}(f_{\beta,p,-}^n(x)) \notin (\mathbf{a}, \mathbf{b}]$ for all $n \ge 0$, i.e., $\tau_{\beta,p,-}(x) \in \Omega^{(\mathbf{a},\mathbf{b}]}$.
As $\tau_{\beta,p,-}(x) = \overline{0} = \tau_{\beta,p,-}(0)$ for all $x < 0$, and $\tau_{\beta,p,-}(x) = \overline{1} = \tau_{\beta,p,-}\big(\frac{\beta}{\beta-1}\big)$ for all $x > \frac{\beta}{\beta-1}$, we obtain that $\tau_{\beta,p,-}(\mathbb{R}) = \tau_{\beta,p,-}\big(\big[0, \frac{\beta}{\beta-1}\big]\big) \subseteq \Omega^{(\mathbf{a},\mathbf{b}]}$.

To show the opposite inclusion, let $\mathbf{c} = c_0 c_1 \cdots \in \Omega^{(\mathbf{a},\mathbf{b}]}$.
If $c_0 = 1$, then let $k_1$ be the length of the maximal common prefix of $\mathbf{c}$ and~$\mathbf{b}$. 
Since $\mathbf{c} > \mathbf{b}$, we have $c_{k_1} = 1$ and $f_{\beta,p,+}^{k_1}(p) < p$. 
Recursively, let $k_{n+1} \ge 1$ be the length of the maximal common prefix of $c_{s_n} c_{s_n+1} \cdots$ and~$\mathbf{b}$, with $s_n = k_1 + \cdots + k_n$.
Then
\begin{align*}
\langle \mathbf{c} \rangle_\beta - p & = \frac{\langle c_{s_1} c_{s_1+1} \cdots \rangle_\beta - f_{\beta,p,+}^{k_1}(p)}{\beta^{s_1}} > \frac{\langle c_{s_1} c_{s_1+1} \cdots \rangle_\beta - p}{\beta^{s_1}} \ge \cdots \\
& \ge \frac{\langle c_{s_n} c_{s_n+1} \cdots \rangle_\beta - f_{\beta,p,+}^{k_n}(p)}{\beta^{s_n}} > \frac{\langle c_{s_n} c_{s_n+1} \cdots \rangle_\beta - p}{\beta^{s_n}}
\end{align*}
for all $n \ge 1$.
Since the latter quantity tends to~$0$ as $n \to \infty$, we have $\langle \mathbf{c} \rangle_\beta > p$.
Similarly, we obtain that $\langle \mathbf{c} \rangle_\beta \le p$ when $c_0 = 0$.
Therefore, we have $\tau_{\beta,p,-}(\langle \mathbf{c} \rangle_\beta) = \mathbf{c}$ for all $\mathbf{c} \in \Omega^{(\mathbf{a},\mathbf{b}]}$, hence $\Omega^{(\mathbf{a},\mathbf{b}]} \subseteq \tau_{\beta,p,-}\big(\big[0, \frac{\beta}{\beta-1}\big]\big)$.

By symmetry, we also get that $\tau_{\beta,p,+}(\mathbb{R}) = \tau_{\beta,p,+}\big(\big[0, \frac{\beta}{\beta-1}\big]\big) = \Omega^{[\mathbf{a},\mathbf{b})}$.

It is well known that $g\big(\tau_{\beta,p,-}\big(\big[0, \frac{\beta}{\beta-1}\big]\big)\big) = \beta = g\big(\tau_{\beta,p,+}\big(\big[0, \frac{\beta}{\beta-1}\big]\big)\big)$; see for example \cite[Proposition~3.7]{Sh}.
Finally,
\[
\Omega^{(\mathbf{a},\mathbf{b})} = \Omega^{(\mathbf{a},\mathbf{b}]} \cup \Omega^{[\mathbf{a},\mathbf{b})}
\]
gives that $g(\Omega^{(\mathbf{a},\mathbf{b})}) = g(\Omega^{(\mathbf{a},\mathbf{b}]}) = g(\Omega^{[\mathbf{a},\mathbf{b})}) = \beta$.
\end{proof}

The address space for negative~$\beta$ can be compared to \cite[Theorem~10]{DK}; see also \cite[Theorem~10]{IS} for the case $p = 1/(\beta-1)$. 

\begin{lemma} \label{l:addneg}
Let $\beta < -1$, $\frac{\beta^2+\beta-1}{\beta^2-1} \le p \le \frac{1}{\beta^2-1}$. 
If $\mathbf{a} := \tau_{\beta,p,-}(p)$ is periodic with odd period length, let $\mathbf{u}$ be its primitive period.  
If $\mathbf{b} := \tau_{\beta,p,+}(p)$ is periodic with odd period length, let $\mathbf{v}$ be its primitive period.  
Then the address space of~$f_{\beta,p,-}$ is
\[
\begin{cases} 
\Omega^{(\mathbf{a},\mathbf{b}]} & \text{if $\mathbf{b}$ is not periodic with odd period length}, \\ \Omega^{(\mathbf{a},\mathbf{b}]} \setminus \{0,1\}^*\, \mathbf{va} & \text{if $\mathbf{b}$ is periodic with odd period length, but $\mathbf{a}$ is not}, \\ \Omega^{(\mathbf{a},\mathbf{b}]} \setminus \{0,1\}^*\, \{\mathbf{v} \overline{\mathbf{u}}, \overline{\mathbf{vu}}\} & \text{if $\mathbf{a}$ and $\mathbf{b}$ are periodic with odd period length},
\end{cases} 
\]
and the address space of~$f_{\beta,p,+}$ is
\[
\begin{cases} 
\Omega^{[\mathbf{a},\mathbf{b})} & \text{if $\mathbf{a}$ is not periodic with odd period length}, \\ 
\Omega^{[\mathbf{a},\mathbf{b})} \setminus \{0,1\}^*\, \mathbf{ub} & \text{if $\mathbf{a}$ is periodic with odd period length, but $\mathbf{b}$ is not}, \\ 
\Omega^{[\mathbf{a},\mathbf{b})} \setminus \{0,1\}^*\, \{\mathbf{u} \overline{\mathbf{v}}, \overline{\mathbf{uv}}\} & \text{if $\mathbf{a}$ and $\mathbf{b}$ are periodic with odd period length},
\end{cases} 
\]
with the alternating lexicographic order on~$\Omega$.
In particular, we have $\mathbf{a} \in \Omega^{(\mathbf{a},\mathbf{b}]}$ and $\mathbf{b} \in \Omega^{[\mathbf{a},\mathbf{b})}$.
Moreover, we have $g\big(\Omega^{(\mathbf{a},\mathbf{b})}\big) = -\beta$. 
\end{lemma}

\begin{proof}
Let $x \in \big[\frac{\beta}{\beta^2-1}, \frac{\beta^2}{\beta^2-1}\big]$, and recall that $f_{\beta,p}\big(\big[\tfrac{\beta}{\beta^2-1}, \tfrac{\beta^2}{\beta^2-1}\big]\big) = \big[\tfrac{\beta}{\beta^2-1}, \tfrac{\beta^2}{\beta^2-1}\big]$ for $\beta < -1$, $\frac{\beta^2+\beta-1}{\beta^2-1} \le p \le \frac{1}{\beta^2-1}$.
Since $f_{\beta,p,-}^n(x) - f_{\beta,p,\pm}^n(p) = \beta^n (x-p)$ when $\tau_{\beta,p,-}(x)$ agrees with $\mathbf{a}$ or~$\mathbf{b}$ on the first $n$ letters, $x < p$ implies that $f_{\beta,p,-}^n(x) \le p < f_{\beta,p,-}^n(p)$ for some even $n \ge 1$ or $f_{\beta,p,-}^n(p) \le p < f_{\beta,p,-}^n(x)$ for some odd $n \ge 1$, hence $\tau_{\beta,p,-}(x) < \mathbf{a}$. 
Similarly, we have $\tau_{\beta,p,-}(x) > \mathbf{b}$ if $x > p$.
This gives that $\tau_{\beta,p,-}(x) \in \Omega^{(\mathbf{a},\mathbf{b}]}$.
As $\tau_{\beta,p,-}(x) = \overline{01} = \tau_{\beta,p,-}\big(\frac{\beta}{\beta^2-1}\big)$ for all $x < \frac{\beta}{\beta^2-1}$, and $\tau_{\beta,p,-}(x) = \overline{10} = \tau_{\beta,p,-}\big(\frac{\beta^2}{\beta-1}\big)$ for all $x > \frac{\beta^2}{\beta^2-1}$, we obtain that $\tau_{\beta,p,-}(\mathbb{R}) = \tau_{\beta,p,-}\big(\big[\frac{\beta}{\beta^2-1}, \frac{\beta^2}{\beta^2-1}\big]) \subseteq \Omega^{(\mathbf{a},\mathbf{b}]}$; in particular, $\mathbf{a} \in \Omega^{(\mathbf{a},\mathbf{b}]}$ and, by symmetry, $\mathbf{b} \in \Omega^{[\mathbf{a},\mathbf{b})}$.
If $\mathbf{b}$ is periodic with odd period length, then $\langle \mathbf{va} \rangle_\beta = p$, but $\mathbf{va}  \ne \mathbf{a} = f_{\beta,p,-}(p)$, hence $\mathbf{va}$ does not occur in the address space of~$f_{\beta,p,-}$. 
If both $\mathbf{a}$ and~$\mathbf{b}$ are periodic with odd period length, then we can also exclude $\overline{\mathbf{vu}}$ because $\langle \overline{\mathbf{vu}} \rangle_\beta = p$.

Now, let $\mathbf{c} = c_0 c_1 \cdots \in \Omega^{(\mathbf{a},\mathbf{b}]}$, $\Omega^{(\mathbf{a},\mathbf{b}]} \setminus \{0,1\}^*\, \mathbf{va}$, and $\Omega^{(\mathbf{a},\mathbf{b}]} \setminus \{0,1\}^*\, \{\mathbf{v} \overline{\mathbf{u}}, \overline{\mathbf{vu}}\}$, respectively.
If $c_0 = 1$, then let $k_1$ be the length of the maximal common prefix of $\mathbf{c}$ and~$\mathbf{b}$. 
If $k_1$ is even, then $c_{k_1} = 1$ and $f_{\beta,p,+}^{k_1}(p) < p$. 
If $k_1$ is odd, then $c_{k_1} = 0$ and $p \le f_{\beta,p,+}^{k_1}(p)$. 
Recursively, let $s_n = k_1 + \cdots + k_n$ and $k_{n+1} \ge 1$ be the length of the maximal common prefix of $c_{s_n} c_{s_n+1} \cdots$ and~$\mathbf{b}$, if $s_n$ is even, the length of the maximal common prefix of $c_{s_n} c_{s_n+1} \cdots$ and~$\mathbf{a}$, if $s_n$ is odd and $c_{s_n} c_{s_n+1} \cdots \ne \mathbf{a}$.
If $c_{s_n} c_{s_n+1} \cdots \ne \mathbf{a}$ all $n \ge 1$, then
\begin{align*}
\langle \mathbf{c} \rangle_\beta - p & = \frac{\langle c_{s_1} c_{s_1+1} \cdots \rangle_\beta - f_{\beta,p,+}^{k_1}(p)}{\beta^{s_1}} \ge \frac{\langle c_{s_1} c_{s_1+1} \cdots \rangle_\beta - p}{\beta^{s_1}} \ge \cdots \\
& \ge \frac{\langle c_{s_n} c_{s_n+1} \cdots \rangle_\beta - f_{\beta,p,\pm}^{k_n}(p)}{\beta^{s_n}} \ge \frac{\langle c_{s_n} c_{s_n+1} \cdots \rangle_\beta - p}{\beta^{s_n}}
\end{align*}
for all $n \ge 2$.
Here, $f_{\beta,p,\pm}^{k_n}(p)$ stands for $f_{\beta,p,+}^{k_n}(p)$ if $s_{n-1}$ is even and for $f_{\beta,p,+}^{k_n}(p)$ if $s_{n-1}$ is odd. 
Since the latter quantity tends to~$0$ as $n \to \infty$, we have $\langle \mathbf{c} \rangle_\beta \ge p$.
This inequality clearly also holds if $c_{s_n} c_{s_n+1} \cdots = \mathbf{a}$ for some $n \ge 1$.
It remains to show that the inequality is strict, i.e., $f_{\beta,p,\pm}^{k_n}(p) \ne p$ for some $n \ge 1$.
If $\mathbf{b}$ is not periodic with odd period length, then this holds for $n = 1$. 
Assume that $\mathbf{b}$ is periodic with primitive period~$\mathbf{v}$ of odd length.
Then $\mathbf{c}$ cannot start with~$\mathbf{vv}$ because this would imply $\mathbf{c} = \overline{\mathbf{v}}$ and thus $\mathbf{c} \notin \Omega^{(\mathbf{a},\mathbf{b}]}$.
Therefore, the only possibility for $f_{\beta,p,+}^{k_1}(p) = p$ is that $k_1 = |\mathbf{v}|$ (and that $c_{k_1} = 0$). 
Since we have excluded that $\mathbf{c} = \mathbf{va}$, we have $c_{k_1} c_{k_1+1} \cdots < \mathbf{a}$. 
If $\mathbf{a}$ is not periodic with odd period length, we have thus $f_{\beta,p,-}^{k_2}(p) \ne p$.
In the remaining case of $\mathbf{a}$ with primitive period~$\mathbf{u}$ of odd length, $c_{k_1} c_{k_1+1} \cdots$ cannot start with~$\mathbf{uu}$, because this would imply that $c_{k_1} c_{k_1+1} \cdots = \mathbf{a}$. 
Thus the only possibility for $f_{\beta,p,-}^{k_2}(p) = p$ is that $k_2 = |\mathbf{u}|$. 
Repeating this argument and since $\mathbf{c} \ne \overline{\mathbf{vu}}$, we obtain that $\langle \mathbf{c} \rangle_\beta > p$.
Since $\langle \mathbf{c} \rangle_\beta \le p$ when $c_0 = 0$, we get that $\tau_{\beta,p,-}(\langle \mathbf{c} \rangle_\beta) = \mathbf{c}$ for all $\mathbf{c} \in \Omega^{(\mathbf{a},\mathbf{b}]}$, $\Omega^{(\mathbf{a},\mathbf{b}]} \setminus \{0,1\}^*\, \mathbf{va}$, and $\Omega^{(\mathbf{a},\mathbf{b}]} \setminus \{0,1\}^*\, \{\mathbf{v} \overline{\mathbf{u}}, \overline{\mathbf{vu}}\}$ respectively, thus this set is equal to $\tau_{\beta,p,-}(\mathbb{R}) = \tau_{\beta,p,-}\big(\big[\frac{\beta}{\beta^2-1}, \frac{\beta^2}{\beta^2-1}\big]\big)$.
By symmetry, $\tau_{\beta,p,+}(\mathbb{R}) = \tau_{\beta,p,+}\big(\big[\frac{\beta}{\beta^2-1}, \frac{\beta^2}{\beta^2-1}\big]\big)$ is $\Omega^{[\mathbf{a},\mathbf{b})}$, $\Omega^{[\mathbf{a},\mathbf{b})} \setminus \{0,1\}^*\, \mathbf{ub}$, and $\Omega^{[\mathbf{a},\mathbf{b})} \setminus \{0,1\}^*\, \{\mathbf{u} \overline{\mathbf{v}}, \overline{\mathbf{uv}}\}$, respectively.

We have $g\big(\tau_{\beta,p,-}\big(\big[\frac{\beta}{\beta^2-1}, \frac{\beta^2}{\beta^2-1}\big]\big)\big) = |\beta| = g\big(\tau_{\beta,p,+}\big(\big[\frac{\beta}{\beta^2-1}, \frac{\beta^2}{\beta^2-1}\big]\big)\big)$ by \cite[Proposition~3.7]{Sh}.
The exponential growth rate of $\Omega^{(\mathbf{a},\mathbf{b}]}$ is the same as that of $\Omega^{(\mathbf{a},\mathbf{b}]} \setminus \{0,1\}^*\, \mathbf{va}$ and $\Omega^{(\mathbf{a},\mathbf{b}]} \setminus \{0,1\}^*\, \{\mathbf{v} \overline{\mathbf{u}}, \overline{\mathbf{vu}}\}$, and a symmetric relation holds for $\Omega^{[\mathbf{a},\mathbf{b})}$.
Together with $\Omega^{(\mathbf{a},\mathbf{b})} = \Omega^{(\mathbf{a},\mathbf{b}]} \cup \Omega^{[\mathbf{a},\mathbf{b})}$, this concludes the proof of the lemma.
\end{proof}

\subsection{Kneading invariant}
The idea for the following lemma goes back to~\cite{MT}; see also \cite{GH,FL}.
Contrary to the cited papers, we work directly with the symbolic space and do not require it to be the address space of some map. 
The main novelty, however, is the treatment of the alternating case.

\begin{lemma} \label{l:K}
Let $\mathbf{a} = a_0 a_1 \cdots$, $\mathbf{b} = b_0 b_1 \cdots \in \Omega$, with $a_0 = 0$, $b_0 = 1$, $\mathbf{a}, \mathbf{b} \in \Omega^{(\mathbf{a},\mathbf{b})}$, and $g(\Omega^{(\mathbf{a},\mathbf{b})}) > 1$, for the lexicographic or alternating lexicographic order on~$\Omega$.
Set
\[
K(z) := \sum_{n=0}^\infty (b_n - a_n)\, z^n = \langle \mathbf{b} \rangle_{1/z} - \langle \mathbf{a} \rangle_{1/z}.
\]
In case of the lexicographic order, $1/g\big(\Omega^{(\mathbf{a},\mathbf{b})}\big)$ is the smallest positive root of~$K(z)$.
In the alternating case, $-1/g\big(\Omega^{(\mathbf{a},\mathbf{b})}\big)$ is the largest negative root of~$K(z)$.
\end{lemma}

\begin{proof}
Define the formal power series 
\[
L(z) := \sum_{n=0}^\infty |L_n| \,  z^n, \quad \text{with} \quad L_n := \big \{ \mathbf{u} \in \{0,1\}^{n}: \, \mathbf u \mathbf c \in \Omega^{(\mathbf a, \mathbf b)} \; \text{for some} \;  \mathbf c \in \Omega \big \},
\]
and
\[
Q(z) := \sum_{n=0}^\infty |Q_j|\,  z^j, \quad \text{with} \quad Q_n := \big\{\mathbf{u} \in \{0,1\}^n:\, \mathbf{ua} \in \Omega^{(\mathbf{a}, \mathbf{b})}\ \text{and}\ \mathbf{ub} \in \Omega^{(\mathbf{a}, \mathbf{b})}\big\}.
\]
We will prove that, for every $n \ge 0$, 
\begin{equation} \label{e:QL}
1 + \sum_{j=0}^n |Q_j| = |L_{n+1}|
\end{equation}
for both lexicographic and alternating lexicographic order, as well as
\begin{equation} \label{e:QL2}
\sum_{j=0}^n (b_{n-j} - a_{n-j})\, |Q_j| = 1 \quad \text{and} \quad \sum_{j=0}^n (-1)^j\, (b_{n-j} - a_{n-j})\, |Q_j| = (-1)^n
\end{equation}
for the lexicographic and alternating lexicographic order, respectively.
By simple formal power series calculations and since $|L_0| = 1$ because $L_0$ contains only the empty word, we obtain from~\eqref{e:QL} that
\[
\frac{1}{1-z}\, \big(1 + Q(z)\big) = \frac{1}{z}\, \big(L(z) - 1\big), 
\]
i.e., $L(z) - \frac{1}{1-z} = Q(z)\, z/(1-z)$.
From~\eqref{e:QL2}, we get that
\[
K(z)\, Q(z) = \frac{1}{1-z} \quad \text{and} \quad K(z)\, Q(-z) = \frac{1}{1+z},
\]
respectively, with the second equation being equivalent to $K(-z)\, Q(z) = \frac{1}{1-z}$.
Then
\begin{equation} \label{e:KLpos}
K(z)\, \bigg(L(z) - \frac{1}{1-z}\bigg) = \frac{z}{(1-z)^2}
\end{equation}
in case of the lexicographic order, and 
\begin{equation} \label{e:KLneg}
K(-z)\, \bigg(L(z) - \frac{1}{1-z}\bigg) = \frac{z}{(1-z)^2}
\end{equation}
in the alternating case. 
Since the radius of convergence of $L(z)$ is~$1/g\big(\Omega^{(\mathbf{a},\mathbf{b})}\big)$ and $L(z)$ has a singularity at~$1/g\big(\Omega^{(\mathbf{a},\mathbf{b})}\big)$, the smallest positive root of $K(z)$ and~$K(-z)$, respectively, is~$1/g\big(\Omega^{(\mathbf{a},\mathbf{b})}\big)$.
It remains to prove~\eqref{e:QL} and~\eqref{e:QL2}.

\smallskip
Consider first the lexicographic order, and order the elements of $L_{n+1}$ lexicographically, from $\mathbf{u}^{(1)} = 0 0 \cdots 0$ to $\mathbf{u}^{(|L_{n+1}|)} = 1 1 \cdots 1$.
For $1 \le k < |L_{n+1}|$, let $\mathbf{u}^{(k)} = u^{(k)}_0 u^{(k)}_1 \cdots u^{(k)}_n$ and
\[ 
\mathbf{v}^{(k)} := u^{(k)}_0 u^{(k)}_1 \cdots u^{(k)}_{j-1}\, b_0 b_1 \cdots b_{n-j},
\]
where $j$ is the minimal integer in $[0,n]$ such that
\begin{equation} \label{e:ukjn}
u^{(k)}_j u^{(k)}_{j+1} \cdots u^{(k)}_n = a_0 a_1 \cdots a_{n-j}.
\end{equation}
(Such an integer exists because $\mathbf{u}^{(k)} \ne 1 1 \cdots 1$ and $\mathbf{u}^{(k)} \mathbf{c} \in \Omega^{(\mathbf{a}, \mathbf{b})}$ for some $\mathbf{c} \in \Omega$.)

We claim that $u^{(k)}_0 u^{(k)}_1 \cdots u^{(k)}_{j-1} \in Q_j$. 
Since $\mathbf{a}, \mathbf{b} \in \Omega^{(\mathbf{a}, \mathbf{b})}$, we have to show that 
\[
u^{(k)}_i u^{(k)}_{i+1} \cdots u^{(k)}_{j-1}\, \mathbf{a} \notin (\mathbf{a}, \mathbf{b}) \quad \text{and} \quad u^{(k)}_i u^{(k)}_{i+1} \cdots u^{(k)}_{j-1}\, \mathbf{b} \notin (\mathbf{a}, \mathbf{b})
\]
for all $0 \le i < j$.
If $u^{(k)}_i = 1$, then \eqref{e:ukjn} and $\mathbf{u}^{(k)} \mathbf{c} \in \Omega^{(\mathbf{a}, \mathbf{b})}$ imply that 
\[
u^{(k)}_i \cdots u^{(k)}_{j-1}\, \mathbf{b} > u^{(k)}_i \cdots u^{(k)}_{j-1}\, \mathbf{a} \ge u^{(k)}_i \cdots u^{(k)}_n\, \mathbf{c} \ge \mathbf{b}.
\]
Assume now that $u^{(k)}_i = 0$.
Then $\mathbf{u}^{(k)} \mathbf{c} \in \Omega^{(\mathbf{a}, \mathbf{b})}$ gives that $u^{(k)}_i \cdots u^{(k)}_n \le a_0 \cdots a_{n-i}$. 
We have $u^{(k)}_i \cdots u^{(k)}_{j-1}\, \mathbf{a} < u^{(k)}_i \cdots u^{(k)}_{j-1}\, \mathbf{b} < \mathbf{a}$ if $u^{(k)}_i \cdots u^{(k)}_{j-1} < a_0 \cdots a_{j-i-1}$, and $u^{(k)}_i \cdots u^{(k)}_n = a_0 \cdots a_{j-i-1}\, a_0 \cdots a_{n-j}$ otherwise, by~\eqref{e:ukjn}.
In the latter case, $a_{j-i} = 0$ implies that $u^{(k)}_i \cdots u^{(k)}_n \ge a_0 \cdots a_{n-i}$, contradicting the minimality of~$j$.
Hence we must have $a_{j-i} = 1$ in this case, thus 
\[
u^{(k)}_i \cdots u^{(k)}_{j-1}\, \mathbf{a} < u^{(k)}_i \cdots u^{(k)}_{j-1}\, \mathbf{b} \le \mathbf{a}.
\]
This proves that $u^{(k)}_0 \cdots u^{(k)}_{j-1} \in Q_j$, thus
\begin{equation} \label{e:LQpos}
L_{n+1} = \{11 \cdots 1\} \cup \bigcup_{j=0}^n \big\{u_0 \cdots u_{j-1}\, a_0 \cdots a_{n-j}:\, u_0 \cdots u_{j-1} \in Q_j\big\}.
\end{equation}
(The inclusion ``$\supseteq$'' is a direct consequence of the definition of $Q_j$ and~$L_{n+1}$.) 
Suppose that $u_0 \cdots u_{i-1}\,  a_0 \cdots a_{n-i} = u_0 \cdots u_{j-1}\, a_0 \cdots a_{n-j}$ with $u_0 \cdots u_{i-1} \in Q_i$, $u_0 \cdots u_{j-1} \in Q_j$, $0 \le i < j \le n$. 
Then we have $a_0 \cdots a_{j-i-1} = u_i \cdots u_{j-1} \in Q_{j-i}$, thus $a_0 \cdots a_{j-i-1}\, \mathbf{b} \le \mathbf{a}$, contradicting that $a_{j-i} = a_0$. 
Therefore, the union in~\eqref{e:LQpos} is disjoint, which gives \eqref{e:QL} for the lexicographic order.

For $1 \le k < |L_{n+1}|$, we have $\mathbf{v}^{(k)} \in L_{n+1}$ since $u^{(k)}_0 \cdots u^{(k)}_{j-1} \in Q_j$.
As there can be no element of $L_{n+1}$ between $\mathbf{u}^{(k)}$ and~$\mathbf{v}^{(k)}$, we obtain that $\mathbf{u}^{(k+1)} = \mathbf{v}^{(k)}$, thus
\[
\sum_{j=0}^n (b_{n-j} - a_{n-j})\, |Q_j| = \sum_{k=1}^{|L_{n+1}|-1} \big(u^{(k+1)}_n - u^{(k)}_n\big) = u^{(|L_{n+1}|)}_n - u^{(1)}_n = 1,
\]
i.e., the left equation in~\eqref{e:QL2} holds for the lexicographic order.

\smallskip
Consider now the alternating lexicographic order, and order the elements of $L_{n+1}$ with respect to this order, from $\mathbf{u}^{(1)} = 0 1 \cdots 0$ to $\mathbf{u}^{(|L_{n+1}|)} = 1 0 \cdots 1$ if $n$ is even, from $\mathbf{u}^{(1)} = 0 1 \cdots 0 1$ to $\mathbf{u}^{(|L_{n+1}|)} = 1 0 \cdots 1 0$ if $n$ is odd.
For $1 \le k < |L_{n+1}|$, let $j$ be the minimal integer in $[0,n]$ such that $u^{(k)}_j \cdots u^{(k)}_n = a_0 \cdots a_{n-j}$ and $j$ is even, or $u^{(k)}_j \cdots u^{(k)}_n = b_0 \cdots b_{n-j}$ and $j$ is odd, with $\mathbf{u}^{(k)} = u^{(k)}_0 \cdots u^{(k)}_n$.
(Such an integer exists because $\mathbf{u}^{(k)}$ is not the maximal element of~$L_{n+1}$.)
Set 
\[ 
\mathbf{v}^{(k)} := u^{(k)}_0 \cdots u^{(k)}_{j-1}\, b_0 \cdots b_{n-j} \quad \text{and} \quad \mathbf{v}^{(k)} := u^{(k)}_0 \cdots u^{(k)}_{j-1}\, a_0 \cdots a_{n-j}
\]
when $j$ is even and odd, respectively.

We claim again that $u^{(k)}_0 \cdots u^{(k)}_{j-1} \in Q_j$. 
Assume w.l.o.g.\ that $j$ is even, the case of odd~$j$ being symmetric.
Let $0 \le i < j$.
If $i$ is even, then we obtain as above that
\[
\begin{array}{cl}
u^{(k)}_i \cdots u^{(k)}_{j-1}\, \mathbf{b} > u^{(k)}_i \cdots u^{(k)}_{j-1}\, \mathbf{a} \ge u^{(k)}_i \cdots u^{(k)}_n\, \mathbf{c} \ge \mathbf{b} & \text{if}\ u^{(k)}_i = 1, \\[1ex]
u^{(k)}_i \cdots u^{(k)}_{j-1}\, \mathbf{a} < u^{(k)}_i \cdots u^{(k)}_{j-1}\, \mathbf{b} \le \mathbf{a} & \text{if}\ u^{(k)}_i = 0.
\end{array}
\]
If $i$ is odd and $u^{(k)}_i = 0$, then
\[
u^{(k)}_i \cdots u^{(k)}_{j-1}\, \mathbf{b} < u^{(k)}_i \cdots u^{(k)}_{j-1}\, \mathbf{a} \le u^{(k)}_i \cdots u^{(k)}_n\, \mathbf{c} \le \mathbf{a}.
\]
Assume now that $i$ is odd and $u^{(k)}_i = 1$.
Then $\mathbf{u}^{(k)} \mathbf{c} \in \Omega^{(\mathbf{a}, \mathbf{b})}$ gives that $u^{(k)}_i \cdots u^{(k)}_n \ge b_0 \cdots b_{n-i}$. 
If $u^{(k)}_i \cdots u^{(k)}_{j-1} = b_0 \cdots b_{j-i-1}$, i.e., $u^{(k)}_i \cdots u^{(k)}_n = b_0 \cdots b_{j-i-1}\, a_0 \cdots a_{n-j}$, then $b_{j-i} = 0$ implies that $u^{(k)}_i \cdots u^{(k)}_n \le b_0 \cdots b_{n-i}$, contradicting the minimality of~$j$.
Hence we must have $b_{j-i} = 1$ in this case, thus 
\[
u^{(k)}_i \cdots u^{(k)}_{j-1}\, \mathbf{a} > u^{(k)}_i \cdots u^{(k)}_{j-1}\, \mathbf{b} \ge \mathbf{b}.
\]
This proves that $u^{(k)}_0 \cdots u^{(k)}_{j-1} \in Q_j$, thus
\begin{multline} \label{e:LQneg}
L_{n+1} = \big\{\mathbf{u}^{(|L_{n+1}|)}\big\} \cup \bigcup_{j=0}^{\lfloor n/2\rfloor} \big\{u_0 \cdots u_{2j-1}\, a_0 \cdots a_{n-2j}:\, u_0 \cdots u_{2j-1} \in Q_{2j}\big\} \\
\cup \bigcup_{j=0}^{\lfloor(n-1)/2\rfloor} \big\{u_0 \cdots u_{2j}\, b_0 \cdots b_{n-2j-1}:\, u_0 \cdots u_{2j} \in Q_{2j+1}\big\}.
\end{multline}
Similarly to the lexicographic case, it is not possible that $u_0 \cdots u_{2i-1}\, a_0 \cdots a_{n-2i} = u_0 \cdots u_{2j-1}\, a_0 \cdots a_{n-2j}$ with $u_0 \cdots u_{2i-1} \in Q_{2i}$, $u_0 \cdots u_{2j-1} \in Q_{2j}$, $0 \le 2i < 2j \le n$. 
If $u_0 \cdots u_{2i-1}\, a_0 \cdots a_{n-2i} = u_0 \cdots u_{2j}\, b_0 \cdots b_{n-2j-1}$ with $u_0 \cdots u_{2i-1} \in Q_{2i}$, $u_0 \cdots u_{2j} \in Q_{2j+1}$, $0 \le 2i \le 2j < n$, then $a_0 \cdots a_{2j-2i} = u_{2i} \cdots u_{2j} \in Q_{2j-2i+1}$, thus $a_0 \cdots a_{2j-2i}\, \mathbf{a} \le \mathbf{a}$, contradicting that $a_{2j-2i+1} = b_0$. 
Other cases of non-empty intersections of two sets on the right hand side of~\eqref{e:LQneg} are excluded symmetrically, thus \eqref{e:QL} holds for the alternating lexicographic order too.

As in the lexicographic case, we have $\mathbf{u}^{(k+1)} = \mathbf{v}^{(k)}$ for $1 \le k < |L_{n+1}|$, thus
\[
\sum_{j=0}^n (-1)^j\, (b_{n-j} - a_{n-j})\, |Q_j| = \sum_{k=1}^{|L_{n+1}|-1} \big(u^{(k+1)}_n - u^{(k)}_n\big) = u^{(|L_{n+1}|)}_n - u^{(1)}_n = (-1)^n,
\]
i.e., the right equation in~\eqref{e:QL2} holds for the alternating lexicographic order.
\end{proof}

\subsection{Monotonicity of $\langle \cdot \rangle_\beta$}

\begin{lemma} \label{l:inc1}
For every $\beta > 2$, $\langle \cdot \rangle_\beta$ and $\langle \cdot \rangle_{-\beta}$ are strictly increasing functions on~$\Omega$, for the lexicographic and alternating lexicographic order, respectively.
\end{lemma}

\begin{proof}
Let $\beta > 2$, and $\mathbf{c}, \mathbf{d} \in \Omega$ with $\mathbf{c} < \mathbf{d}$ (for the lexicographic or alternating lexicographic order). 
By removing the maximum initial portion of the strings where $\mathbf{c}$ and
$\mathbf{d}$ are equal, and exchanging the role of~$\mathbf{c}$ and~$\mathbf{d}$ if the length of this portion is odd in the alternating case, we may assume w.l.o.g.\ that $\mathbf{c}$ starts with~$0$ and $\mathbf{d}$ starts with~$1$.
Then we have 
\[
\langle \mathbf{c} \rangle_\beta \le \langle 0\, \overline{1} \rangle_\beta = \frac{1}{\beta-1} < 1 = \langle 1\, \overline{0} \rangle_\beta \le \langle \mathbf{d} \rangle_\beta
\]
in the lexicographic case, and
\[
\langle \mathbf{c} \rangle_{-\beta} \le \langle 0\, \overline{01} \rangle_{-\beta} = \frac{1}{\beta^2-1} < \frac{\beta^2-\beta-1}{\beta^2-1} = \langle 1\, \overline{10} \rangle_{-\beta} \le \langle \mathbf{d} \rangle_{-\beta}
\]
in the alternating case.
This proves the lemma.
\end{proof}

\begin{lemma} \label{l:inc2}
Let $\mathbf{a} \in 0\, \Omega$, $\mathbf{b} \in 1\, \Omega$, with $\mathbf{a}, \mathbf{b} \in \Omega^{(\mathbf{a},\mathbf{b})}$ and $g\big(\Omega^{(\mathbf{a},\mathbf{b})}\big) > 1$, where $\Omega$ is equipped with the lexicographic or the alternating lexicographic order.

In case of the lexicographic order, $\langle \cdot \rangle_{g(\Omega^{(\mathbf{a},\mathbf{b})})}$ is an increasing function on~$\Omega^{(\mathbf{a},\mathbf{b})}$.

In the alternating case, $\langle \cdot \rangle_{-g(\Omega^{(\mathbf{a},\mathbf{b})})}$ is an increasing function on~$\Omega^{(\mathbf{a},\mathbf{b})}$.
\end{lemma}

\begin{proof}
In the following, we assume that $\Omega$ is equipped with the lexicographic order; for the alternating case, we only have to change $\langle \cdot \rangle_\beta$ to $\langle \cdot \rangle_{-\beta}$.

We show that $\langle \cdot \rangle_\beta$ is strictly increasing on $\Omega^{(\mathbf{a},\mathbf{b})}$ for all $\beta > g\big(\Omega^{(\mathbf{a},\mathbf{b})}\big)$.
By Lemma~\ref{l:inc1}, this is true for $\beta > 2$. 
From Lemma~\ref{l:K}, we know that $\langle \mathbf{a} \rangle_\beta \ne \langle \mathbf{b} \rangle_\beta$ for all $\beta > g\big(\Omega^{(\mathbf{a},\mathbf{b})}\big)$.
By the continuity of $\langle \mathbf{a} \rangle_\beta$ and $\langle \mathbf{b} \rangle_\beta$ as functions of~$\beta > 1$ and since $\langle \mathbf{a} \rangle_\beta < \mathbf{b} \rangle_\beta$ for all $\beta > 2$, we obtain that $\langle \mathbf{a} \rangle_\beta < \langle \mathbf{b} \rangle_\beta$ for all $\beta > g\big(\Omega^{(\mathbf{a},\mathbf{b})}\big)$.

Assume that $\mathbf{c} < \mathbf{d}$, but $\langle \mathbf{c} \rangle_\beta \ge \langle \mathbf{d} \rangle_\beta$ for some $\mathbf{c}, \mathbf{d} \in \Omega^{(\mathbf{a},\mathbf{b})}$, $\beta > g\big(\Omega^{(\mathbf{a},\mathbf{b})}\big)$.
By removing the longest common prefix of $\mathbf{c}$ and~$\mathbf{d}$, we may assume w.l.o.g.\ that $\mathbf{c}$ starts with~$0$ and $\mathbf{d}$ starts with~$1$.
(In the alternating case, we also have to exchange the role of~$\mathbf{c}$ and~$\mathbf{d}$ if the length of this prefix is odd.)
Consider
\[
\Delta(\beta) := \max \big\{\langle \mathbf{c} \rangle_\beta - \langle \mathbf{d} \rangle_\beta:\, \mathbf{c}, \mathbf{d} \in \Omega^{(\mathbf{a},\mathbf{b})},\, \mathbf{c} \le \mathbf{a},\, \mathbf{b} \le \mathbf{d}\big\} \qquad (\beta > 1).
\]
This function is well defined because $\{\langle \mathbf{c} \rangle_\beta:\, \mathbf{c} \in \Omega^{(\mathbf{a},\mathbf{b})}\}$ is compact; $\Delta(\cdot)$~is also continuous.
By Lemma~\ref{l:inc1}, we have $\Delta(\beta) < 0$ for $\beta > 2$, and we have assumed that $\Delta(\beta) \ge 0$ for some $\beta > g\big(\Omega^{(\mathbf{a},\mathbf{b})}\big)$. 
Therefore, there exists $\beta > g\big(\Omega^{(\mathbf{a},\mathbf{b})}\big)$ such that $\Delta(\beta) = 0$.
Fix this~$\beta$, and choose $\mathbf{c}, \mathbf{d} \in \Omega^{(\mathbf{a},\mathbf{b})}$ with $\mathbf{c} \le \mathbf{a}$,  $\mathbf{b} \le \mathbf{d}$, and $\langle \mathbf{c} \rangle_\beta = \langle \mathbf{d} \rangle_\beta$.
Since $\langle \mathbf{a} \rangle_\beta < \langle \mathbf{b} \rangle_\beta$, we have $\langle \mathbf{c} \rangle_\beta > \langle \mathbf{a} \rangle_\beta$ or $\langle \mathbf{b} \rangle_\beta > \langle \mathbf{d} \rangle_\beta$.
If $\langle \mathbf{c} \rangle_\beta > \langle \mathbf{a} \rangle_\beta$, then removing the longest common prefix of $\mathbf{a}$ and~$\mathbf{c}$ gives sequences $\mathbf{c}', \mathbf{d}' \in \Omega^{(\mathbf{a},\mathbf{b})}$ with $\mathbf{c}' \le \mathbf{a}$, $\mathbf{d'} \ge \mathbf{b}$, and $\langle \mathbf{c}' \rangle_\beta > \langle \mathbf{d}' \rangle_\beta$, contradicting that $\Delta(\beta) = 0$. 
Similarly, $\langle \mathbf{b} \rangle_\beta > \langle \mathbf{d} \rangle_\beta$ leads to a contradiction.
Therefore, we have shown that $\langle \cdot \rangle_\beta$ is strictly increasing on~$\Omega^{(\mathbf{a},\mathbf{b})}$ for all $\beta > g\big(\Omega^{(\mathbf{a},\mathbf{b})}\big)$.

By continuity in~$\beta$, we obtain that $\langle \cdot \rangle_{g(\Omega^{(\mathbf{a},\mathbf{b})})}$ is increasing on~$\Omega^{(\mathbf{a},\mathbf{b})}$.
\end{proof}

\subsection{Periodic critical itineraries}

Next we show that condition (\ref{i:adm3}) of Definition~\ref{def:admissible} is violated when  (\ref{i:adm1}) and  (\ref{i:adm2}) hold but $\mathbf{a}$ and $\mathbf{b}$ are not the critical itineraries of some~$f_{\beta,p}$.

\begin{lemma} \label{l:notcritical}
Let $\mathbf{a} \in 0\, \Omega$, $\mathbf{b} \in 1\, \Omega$, with $\mathbf{a} \in \Omega^{(\mathbf{a},\mathbf{b}]}$, $\mathbf{b} \in \Omega^{[\mathbf{a},\mathbf{b})}$, and $g\big(\Omega^{(\mathbf{a},\mathbf{b})}\big) > 1$.
Set $\beta := g\big(\Omega^{(\mathbf{a},\mathbf{b})}\big)$ when $\Omega$ is equipped with the lexicographic order, $\beta := -g\big(\Omega^{(\mathbf{a},\mathbf{b})}\big)$ in case of the alternating lexicographic order, and $p := \langle \mathbf{a} \rangle_\beta$. 
If $\mathbf{a} \ne \tau_{\beta,p,-}(p)$ or $\mathbf{b} \ne \tau_{\beta,p,+}(p)$, then we have $\mathbf{a}, \mathbf{b} \in \{\mathbf{u}, \mathbf{v}\}^\omega$ for some $\mathbf{u} \in 0\, \{0,1\}^*$, $\mathbf{v} \in 1\, \{0,1\}^*$, with $\overline{\mathbf{u}} \in \Omega^{(\overline{\mathbf{u}},\overline{\mathbf{v}}]}$, $\overline{\mathbf{v}} \in \Omega^{[\overline{\mathbf{u}},\overline{\mathbf{v}})}$, $g(\Omega^{(\overline{\mathbf{u}},\overline{\mathbf{v}})}) = g(\Omega^{(\mathbf{a},\mathbf{b})})$, and $\mathbf{a} \ne \overline{\mathbf{u}}$ or $\mathbf{b} \ne \overline{\mathbf{v}}$.
\end{lemma}

\begin{proof}
First note that $\langle \mathbf{b} \rangle_\beta = \langle \mathbf{a} \rangle_\beta = p$ by Lemma~\ref{l:K}.
Assume that $\mathbf{a} \ne \tau_{\beta,p,-}(p)$; the case of $\mathbf{b} \ne \tau_{\beta,p,+}(p)$ is symmetric.
Write $\mathbf{a} = a_0 a_1 \cdots$ and $\tau_{\beta,p,-}(p) = c_0 c_1 \cdots$.
If $a_0 \cdots a_{n-1} = c_0 \cdots c_{n-1}$, $n \ge 0$, then $f_{\beta,p,-}^n(p) = \langle c_n c_{n+1} \cdots \rangle_\beta = \langle a_n a_{n+1} \cdots \rangle_\beta$.
By Lemma~\ref{l:inc2}, $\langle a_n a_{n+1} \cdots \rangle_\beta < p$ implies that $a_n = 0$, and $\langle a_n a_{n+1} \cdots \rangle_\beta > p$ implies that $a_n = 1$; hence $a_n \neq c_n$ is possible only when $f_{\beta,p,-}^n(p) = p$. 
Since $\mathbf{a} \ne \tau_{\beta,p,-}(p)$ and $a_0 = 0 = c_0$, the latter case occurs for some $n \ge 1$.
Choose $m \ge 1$ minimal such that $f_{\beta,p,-}^m(p) = p$, and set $\mathbf{u} := a_0 \cdots a_{m-1}$; then $\tau_{\beta,p,-}(p) = \overline{\mathbf{u}}$.
Since $\mathbf{a} \ne \tau_{\beta,p,-}(p)$, there exists some $\ell \ge 1$ such that $a_1 \cdots a_{\ell m-1} = \mathbf{u} \cdots \mathbf{u} = c_0 \cdots c_{\ell m-1}$ and $a_{\ell m} = 1$.
Then we have $\langle a_{\ell m} a_{\ell m+1} \cdots \rangle_\beta = p = \langle \mathbf{b} \rangle_\beta$ and $a_{\ell m} a_{\ell m+1} \cdots > \mathbf{b}$ since $\mathbf{a} \in \Omega^{(\mathbf{a},\mathbf{b}]}$.
Similarly as for $a_0 a_1 \cdots$ and $c_0 c_1 \cdots$, we obtain that there exists some $n \ge 1$ such that $a_{\ell m} \cdots a_{\ell m+n-1} = b_0 \cdots b_{n-1}$ and $\langle b_n b_{n+1} \cdots \rangle_\beta = p$, hence we also have $f_{\beta,p,+}^k(p) = p$ for some $k \ge 1$. 
Let $j \ge 1$ be minimal such that $f_{\beta,p,+}^j(p) = p$, and set $\mathbf{v} := b_0 \cdots b_{j-1}$. 
Then we have $\tau_{\beta,p,+}(p) = \overline{\mathbf{v}}$ and $\mathbf{a}, \mathbf{b} \in \{\mathbf{u}, \mathbf{v}\}^\omega$.

Since $\overline{\mathbf{u}}, \overline{\mathbf{v}}$ are the itineraries of $f_{\beta,p,\pm}(p)$, we have $\overline{\mathbf{u}} \in \Omega^{(\overline{\mathbf{u}},\overline{\mathbf{v}}]}$, $\overline{\mathbf{v}} \in \Omega^{[\overline{\mathbf{u}},\overline{\mathbf{v}})}$, and $g\big(\Omega^{(\overline{\mathbf{u}},\overline{\mathbf{v}})}\big) = |\beta| = g\big(\Omega^{(\mathbf{a},\mathbf{b})}\big)$ by Lemma~\ref{l:addpos} and~\ref{l:addneg}, respectively.
\end{proof}

\begin{lemma} \label{l:val}
Let $\mathbf{u} \in 0\, \{0,1\}^*$, $\mathbf{v} \in 1\, \{0,1\}^*$, $|\beta| > 1$, with $\langle \overline{\mathbf{u}} \rangle_\beta = \langle \overline{\mathbf{v}} \rangle_\beta =: p$. 

Then $\langle \mathbf{c} \rangle_\beta = p$ for all $\mathbf{c} \in \{\mathbf{u}, \mathbf{v}\}^\omega$.

If $\tau_{\beta,p,-}(p) \in \{\mathbf{u}, \mathbf{v}\}^\omega$, then $\tau_{\beta,p,-}(p) = \overline{\mathbf{u}}$. 

If $\tau_{\beta,p,+}(p) \in \{\mathbf{u}, \mathbf{v}\}^\omega$, then $\tau_{\beta,p,+}(p) = \overline{\mathbf{v}}$. 
\end{lemma}

\begin{proof}
Let $\mathbf{c} = c_0 c_1 \cdots \in \Omega$ with $\mathbf{c} \in \{\mathbf{u}, \mathbf{v}\}^\omega$. 
If $\mathbf{c}$ starts with~$\mathbf{u}$, then 
\[
\langle \mathbf{c} \rangle_\beta - p = \langle \mathbf{c} \rangle_\beta - \langle \overline{\mathbf{u}} \rangle_\beta = \frac{\langle c_{|\mathbf{u}|} c_{|\mathbf{u}|+1} \cdots \rangle_\beta - \langle \overline{\mathbf{u}} \rangle_\beta}{\beta^{|\mathbf{u}|}} = \frac{\langle c_{|\mathbf{u}|} c_{|\mathbf{u}|+1} \cdots \rangle_\beta - p}{\beta^{|\mathbf{u}|}}.
\]
Similarly, we have $\langle \mathbf{c} \rangle_\beta - p = (\langle c_{|\mathbf{v}|} c_{|\mathbf{v}|+1} \cdots \rangle_\beta - p)/\beta^{|\mathbf{v}|}$ if $\mathbf{c}$ starts with~$\mathbf{v}$, and thus 
\[
\langle \mathbf{c} \rangle_\beta - p = \frac{\langle c_n c_{n+1} \cdots \rangle_\beta - p}{\beta^n}
\]
for each $n \ge 1$ such that $c_0 \cdots c_{n-1} \in \{\mathbf{u}, \mathbf{v}\}^*$. 
Since $n$ is unbounded and $\langle c_n c_{n+1} \cdots \rangle_\beta$ is bounded, we get that~$\langle \mathbf{c} \rangle_\beta = p$.

If $\tau_{\beta,p,-}(p) \in \{\mathbf{u}, \mathbf{v}\}^\omega$, then $\tau_{\beta,p,-}(p)$ starts with~$\mathbf{u}$, thus $f_{\beta,p,-}^{|\mathbf{u}|}(p) = p$ and $\tau_{\beta,p,-}(p) = \overline{\mathbf{u}}$.
Similarly, $\tau_{\beta,p,+}(p) \in \{\mathbf{u}, \mathbf{v}\}^\omega$ implies that $\mathbf{b} = \overline{\mathbf{v}}$. 
\end{proof}

\subsection{Proof of Theorem~\ref{t:main}}
Assume first that $\mathbf{a} = \tau_{\beta,p,-}(p)$ and $\mathbf{b} = \tau_{\beta,p,+}(p)$ for some $\beta, p$. 
Then condition (\ref{i:adm1}) of Definition~\ref{def:admissible} holds by Lemmas~\ref{l:addpos} and~\ref{l:addneg}, respectively.
These lemmas also give that $g\big(\Omega^{(\mathbf{a},\mathbf{b})}\big) = |\beta| > 1$.
If $\mathbf{u}, \mathbf{v}$ are as in condition~(\ref{i:adm3}) of Definition~\ref{def:admissible}, then $\langle \overline{\mathbf{u}} \rangle_\beta = \langle \overline{\mathbf{v}} \rangle_\beta$ by Lemma~\ref{l:K}, and this value equals~$p$ by Lemma~\ref{l:val}.
Lemma~\ref{l:val} also gives that $\mathbf{a} = \overline{\mathbf{u}}$ and $\mathbf{b} = \overline{\mathbf{v}}$, thus condition~(\ref{i:adm3}) of Definition~\ref{def:admissible} holds.
Therefore, the pair $(\mathbf{a}; \mathbf{b})$ is lex-admissible if $\beta > 1$ and alt-admissible if $\beta < -1$.

Now, let $(\mathbf{a}; \mathbf{b})$ be lex-admissible and $\beta := g(\Omega^{(\mathbf{a},\mathbf{b})})$, or alt-admissible and $\beta := -g(\Omega^{(\mathbf{a},\mathbf{b})})$.
By Lemma~\ref{l:notcritical} and condition~(\ref{i:adm3}) of Definition~\ref{def:admissible}, we have $\mathbf{a} = \tau_{\beta,p,-}(p)$ and $\mathbf{b} = \tau_{\beta,p,+}(p)$, with 
$p := \langle \mathbf{a} \rangle_\beta$.
Moreover, we have $1 \le \langle \mathbf{b} \rangle_\beta = p = \langle \mathbf{a} \rangle_\beta \le \frac{1}{\beta-1}$ in case $\beta > 1$, and $1 + \frac{\beta}{\beta^2-1} \le \langle \mathbf{b} \rangle_\beta = p = \langle \mathbf{a} \rangle_\beta \le \frac{1}{\beta^2-1}$ in case $\beta < -1$.

\subsection{Proof of Corollary~\ref{c:main}}
Let $\beta > 1$, $0 \le \alpha \le 2 - \beta$, set $p := \frac{\alpha}{\beta-1} + 1$ and define the map $\varphi:\, [0,1] \to [\beta (p-1), \beta p]$, $x \mapsto \beta (x+p-1)$.
Then $T_{\beta,\alpha,\pm} = \varphi^{-1} \circ f_{\beta,p,\pm} \circ \varphi$.
Therefore, $\tau_{\beta,p,-}(p)$ and $\tau_{\beta,p,+}(p)$ are the itineraries of $\varphi^{-1}(p) = (1-\alpha)/\beta$ under $T_{\beta,\alpha,-}$ and~$T_{\beta,\alpha,+}$.
Since each $p \in \big[1, \frac{1}{\beta-1}\big]$ can be written as $p = \frac{\alpha}{\beta-1} + 1$ with $0 \le \alpha \le 2 - \beta$, the critical itineraries of $f_{\beta,p}$ are exactly those of $\beta x + \alpha \bmod 1$. 

Now, let $\beta < -1$, $-\beta - 1 < \alpha < 1$, set $p := \frac{\alpha}{1-\beta}$ and define the map $\varphi:\, [0,1] \to [\beta p, \beta (p-1)]$, $x \mapsto \beta (p-x)$.
Then $T_{\beta,\alpha,\pm} = \varphi^{-1} \circ f_{\beta,p,\pm} \circ \varphi$.
Therefore, the itineraries of $\varphi^{-1}(p) = -\alpha/\beta$ under $T_{\beta,\alpha,-}$ and~$T_{\beta,\alpha,+}$ are $\mathbf{a} = \tau_{\beta,p,-}(p)$ and $\mathbf{b} = \tau_{\beta,p,+}(p)$.
Hence $(\mathbf{a}; \mathbf{b})$ is alt-admissible.
Moreover, we have $\mathbf{a} \in 00\, \Omega$, $\mathbf{b} \in 11\, \Omega$, and $\frac{\beta+1}{\beta-1} < p < \frac{1}{1-\beta}$ holds if and only if $f_{\beta,p,-}^2(p) = \beta^2 p < \beta (p-1) = f_{\beta,p,+}(p)$ and $f_{\beta,p,+}^2(p) = \beta^2 p - \beta^2 - \beta > \beta p = f_{\beta,p,-}(p)$, i.e., $S^2(\mathbf{a}) < S(\mathbf{b})$ and $S^2(\mathbf{b}) > S(\mathbf{a})$.

On the other hand, let $(\mathbf{a}; \mathbf{b})$ be alt-admissible, $S^2(\mathbf{a}) < S(\mathbf{b})$ and $S^2(\mathbf{b}) > S(\mathbf{a})$, and set $\beta := -g\big(\Omega^{(\mathbf{a},\mathbf{b})}\big)$, $\alpha := p\, (1-\beta)$ with $p := \langle \mathbf{a} \rangle_\beta$. 
We first show that $\mathbf{a} \in 00\, \Omega$ and $\mathbf{b} \in 11\, \Omega$.
Assume that $\mathbf{a}$ starts with~$01$. 
Then each $0$ in $\mathbf{a}$ or~$\mathbf{b}$ is followed by a~$1$. 
Moreover, $S^2(\mathbf{b}) > S(\mathbf{a})$ implies that $\mathbf{b}$ starts with $111$ or~$101$. 
In the latter case we have $\mathbf{a} = \overline{01}$ and $\mathbf{b} = \overline{10}$, contradicting that $g\big(\Omega^{(\mathbf{a},\mathbf{b})}\big) > 1$.
If $\mathbf{b}$ starts with $111$, then we must have $\mathbf{b} = \overline{1}$, hence each $1$ in~$\mathbf{a}$ is followed by a~$0$ as otherwise $\mathbf{a}$ would end with~$\mathbf{b}$.
This gives that $\mathbf{a} = \overline{01}$, also contradicting that $g\big(\Omega^{(\mathbf{a},\mathbf{b})}\big) > 1$.
By symmetry, $\mathbf{b}$~cannot start with~$10$.
Now, we obtain from the considerations in the previous paragraph that $\mathbf{a}, \mathbf{b}$ are the critical itineraries of $\beta x + \alpha \bmod 1$, with $\beta < -1$, $-\beta - 1 < \alpha < 1$.

\section{Exponential growth rates}
\label{sec:expon-growth-rates}

\subsection{Lexicographic order}
The following lemma is similar to \cite[Theorem~3]{SS}.

\begin{lemma} \label{l:entpos}
Let $\mathbf{a} \in 0\, \Omega$, $\mathbf{b} \in 1\, \Omega$, $\mathbf{u} \in 01\, \{0,1\}^*$, $\mathbf{v} \in 10\, \{0,1\}^*$, such that $\mathbf{a}, \mathbf{b} \in \{\mathbf{u}, \mathbf{v}\}^\omega$, $\mathbf{a} \in \Omega^{(\mathbf{a},\mathbf{b}]}$, $\mathbf{b} \in \Omega^{[\mathbf{a},\mathbf{b})}$, $\overline{\mathbf{u}} \in \Omega^{(\overline{\mathbf{u}},\overline{\mathbf{v}}]}$, and $\overline{\mathbf{v}} \in \Omega^{[\overline{\mathbf{u}},\overline{\mathbf{v}})}$, where $\Omega$ is equipped with the lexicographic order.
Then 
\[
g\big(\Omega^{(\overline{\mathbf{u}},\overline{\mathbf{v}})}\big)  = g\big(\Omega^{(\mathbf{a},\mathbf{b})}\big).
\]
\end{lemma}

\begin{proof}
Since
\[
\Omega^{(\overline{\mathbf{u}},\overline{\mathbf{v}})} \subseteq \Omega^{(\mathbf{a},\mathbf{b})} \subseteq \Omega^{(\mathbf{u}\overline{\mathbf{v}},\mathbf{v}\overline{\mathbf{u}})},
\]
it suffices to show that $g\big(\Omega^{(\mathbf{u}\overline{\mathbf{v}},\mathbf{v}\overline{\mathbf{u}})}\big) = g\big(\Omega^{(\overline{\mathbf{u}},\overline{\mathbf{v}})}\big)$.

First we prove that  
\begin{equation} \label{e:capuvpos}
\Omega^{(\mathbf{u}\overline{\mathbf{v}},\mathbf{v}\overline{\mathbf{u}})} \cap [\overline{\mathbf{u}}, \overline{\mathbf{v}}] \subseteq \{\mathbf{u}, \mathbf{v}\}^\omega.
\end{equation}
Let $\mathbf{c} \in \Omega^{(\mathbf{u}\overline{\mathbf{v}},\mathbf{v}\overline{\mathbf{u}})} \cap [\overline{\mathbf{u}}, \overline{\mathbf{v}}]$, and consider an arbitrary decomposition $\mathbf{c} = \mathbf{wd}$ with $\mathbf{w} \in \{\mathbf{u}, \mathbf{v}\}^*$, $\mathbf{d} \in \Omega$.
Assume that $\mathbf{d}$ starts with~$0$.
Then $\mathbf{c} \in \Omega^{(\mathbf{u}\overline{\mathbf{v}},\mathbf{v}\overline{\mathbf{u}})}$ implies that $\mathbf{d} \le \mathbf{u} \overline{\mathbf{v}}$ and that $\mathbf{d} \ge \overline{\mathbf{u}}$ when $\mathbf{w}$ ends with a word in $\mathbf{v}\, \{\mathbf{u}\}^*$. 
If $\mathbf{w}$ contains no occurrence of~$\mathbf{v}$, i.e., $\mathbf{w} \in \{\mathbf{u}\}^*$, then we get $\mathbf{d} \ge \overline{\mathbf{u}}$ from $\mathbf{c} \in [\overline{\mathbf{u}}, \overline{\mathbf{v}}]$.
Therefore, we have $\overline{\mathbf{u}} \le \mathbf{d} \le \mathbf{u} \overline{\mathbf{v}}$, hence $\mathbf{d}$ starts with~$\mathbf{u}$. 
Symmetrically, we obtain that  $\mathbf{d}$ starts with~$\mathbf{v}$ whenever it starts with~$1$. 
Therefore, \eqref{e:capuvpos} holds, and
\[
\Omega^{(\overline{\mathbf{u}},\overline{\mathbf{v}})} \subseteq \Omega^{(\mathbf{u}\overline{\mathbf{v}},\mathbf{v}\overline{\mathbf{u}})} \subseteq \Omega^{(\overline{\mathbf{u}},\overline{\mathbf{v}})} \cup \bigcup_{n=0}^\infty \Omega^{(\overline{\mathbf{u}},\overline{\mathbf{v}})}_n\, \{\mathbf{u}, \mathbf{v}\}^\omega,
\]
where $\Omega^{(\overline{\mathbf{u}},\overline{\mathbf{v}})}_n$ denotes the set of length $n$ prefixes of words in~$\Omega^{(\overline{\mathbf{u}},\overline{\mathbf{v}})}$.
This implies that $g\big(\Omega^{(\overline{\mathbf{u}},\overline{\mathbf{v}})}\big) \le g\big(\Omega^{(\mathbf{u}\overline{\mathbf{v}},\mathbf{v}\overline{\mathbf{u}})}\big) \le \max\big( g\big(\Omega^{(\overline{\mathbf{u}},\overline{\mathbf{v}})}\big), g\big(\bigcup_{n=0}^\infty S^n \{\mathbf{u}, \mathbf{v}\}^\omega\big) \big)$. 
Hence we only have to show that $g\big(\bigcup_{n=0}^\infty S^n \{\mathbf{u}, \mathbf{v}\}^\omega\big) \le g\big(\Omega^{(\overline{\mathbf{u}},\overline{\mathbf{v}})}\big)$.

Next we prove that
\begin{equation} \label{e:incrs}
\{\mathbf{r}, \mathbf{s}\}^\omega \subseteq \Omega^{(\overline{\mathbf{u}},\overline{\mathbf{v}})},
\end{equation}
where $\mathbf{r}$ is the longest common prefix of $\overline{\mathbf{u}}$ and~$S^n(\overline{\mathbf{v}})$, with $n \ge 1$ such that $S^n(\overline{\mathbf{v}})$ is maximal among all suffixes of~$\mathbf{v}$ starting with~$0$, and $\mathbf{s}$ is the longest common prefix of $\overline{\mathbf{v}}$ and~$S^m(\overline{\mathbf{u}})$, with $m \ge 1$ such that $S^m(\overline{\mathbf{u}})$ is minimal among all suffixes of~$\mathbf{u}$ starting with~$1$.
Note that $\mathbf{r}$ and~$\mathbf{s}$ are finite words because $S^n(\overline{\mathbf{v}}) < \overline{\mathbf{u}}$ by $\overline{\mathbf{v}} \in \Omega^{[\overline{\mathbf{u}},\overline{\mathbf{v}})}$ and $S^m(\overline{\mathbf{u}}) > \overline{\mathbf{v}}$ by $\overline{\mathbf{u}} \in \Omega^{(\overline{\mathbf{u}},\overline{\mathbf{v}}]}$. 
Moreover, $S^n(\overline{\mathbf{v}})$~starts with~$\mathbf{r} 0$ and $\overline{\mathbf{u}}$ starts with~$\mathbf{r} 1$, $\overline{\mathbf{v}}$~starts with~$\mathbf{s} 0$ and $S^m(\overline{\mathbf{u}})$~starts with~$\mathbf{s} 1$.
Therefore, we have 
\[
S^n(\overline{\mathbf{v}}) \le \mathbf{r} S^n(\overline{\mathbf{v}}) \le \cdots \le \overline{\mathbf{r}} < \overline{\mathbf{s}} \le \cdots \le \mathbf{s} S^m(\overline{\mathbf{u}}) \le S^m(\overline{\mathbf{u}}).
\] 
Since $S^n(\overline{\mathbf{v}}) \in \Omega^{(\overline{\mathbf{u}},\overline{\mathbf{v}})}$ and $\overline{\mathbf{u}} \in \Omega^{(\overline{\mathbf{u}},\overline{\mathbf{v}})}$, we have $\mathbf{rc} \in \Omega^{(\overline{\mathbf{u}},\overline{\mathbf{v}})}$ for all $\mathbf{c} \in \Omega^{(\overline{\mathbf{u}},\overline{\mathbf{v}})}$ with $\mathbf{c} \in [S^{n+|\mathbf{r}|}(\overline{\mathbf{v}}), S^{|\mathbf{r}|}(\overline{\mathbf{u}})] \supseteq [S^n(\overline{\mathbf{v}}), S^m(\overline{\mathbf{u}})]$. 
Symmetrically, we also have $\mathbf{sc} \in \Omega^{(\overline{\mathbf{u}},\overline{\mathbf{v}})}$ for all $\mathbf{c} \in \Omega^{(\overline{\mathbf{u}},\overline{\mathbf{v}})} \cap [S^n(\overline{\mathbf{v}}), S^m(\overline{\mathbf{u}})]$, thus 
\[
\mathbf{rc}, \mathbf{sc} \in \Omega^{(\overline{\mathbf{u}},\overline{\mathbf{v}})} \cap [S^n(\overline{\mathbf{v}}), S^m(\overline{\mathbf{u}})] \quad \text{for all} \quad  \mathbf{c} \in \Omega^{(\overline{\mathbf{u}},\overline{\mathbf{v}})} \cap [S^n(\overline{\mathbf{v}}), S^m(\overline{\mathbf{u}})].
\]
By compactness of~$\Omega$, we get that~\eqref{e:incrs} holds, thus $g\big(\bigcup_{n=0}^\infty S^n \{\mathbf{r}, \mathbf{s}\}^\omega\big) \le g\big(\Omega^{(\overline{\mathbf{u}},\overline{\mathbf{v}})}\big)$.

It remains to prove that $g\big(\bigcup_{n=0}^\infty S^n \{\mathbf{u}, \mathbf{v}\}^\omega\big) \le g\big(\bigcup_{n=0}^\infty S^n \{\mathbf{r}, \mathbf{s}\}^\omega\big)$.
To this end, we show that $\min(|\mathbf{r}|, |\mathbf{s}|) \le \min(|\mathbf{u}|, |\mathbf{v}|)$ and $\max(|\mathbf{r}|, |\mathbf{s}|) \le \max(|\mathbf{u}|, |\mathbf{v}|)$. 
The latter inequality follows from the definition of $\mathbf{r}$ and~$\mathbf{s}$. 
If $|\mathbf{u}| = |\mathbf{v}|$, then the former inequality also holds.
It remains to consider the case $|\mathbf{u}| \ne |\mathbf{v}|$; w.l.o.g.\ $|\mathbf{u}| < |\mathbf{v}|$.
Suppose that $|\mathbf{r}| > |\mathbf{u}|$. 
Then, for some $k \ge n$, $S^k(\overline{\mathbf{v}})$ starts with~$\mathbf{ut} 0$, with some prefix~$\mathbf{t}$ of~$\mathbf{u}$, while $\overline{\mathbf{u}}$ starts with~$\mathbf{ut} 1$.
Write $S^m(\overline{\mathbf{u}}) = \mathbf{w} 1 \mathbf{c}$, with $\mathbf{w}$ a suffix of $\mathbf{ut}$ of length at most~$|\mathbf{u}|$.
As $S^{k+|\mathbf{ut}|-|\mathbf{w}|}(\overline{\mathbf{v}})$ starts with~$\mathbf{w} 0$, and $\overline{\mathbf{v}} \le S^{k+|\mathbf{ut}|-|\mathbf{w}|}(\overline{\mathbf{v}})$, we obtain that $|\mathbf{s}| \le |\mathbf{w}| \le |\mathbf{u}|$. 
Therefore, we have $\min(|\mathbf{r}|, |\mathbf{s}|) \le \min(|\mathbf{u}|, |\mathbf{v}|)$, which concludes the proof of the lemma.
\end{proof}

\subsection{Proof of Theorem~\ref{t:pos}}
We show first that (\ref{i:pos3}) implies (\ref{i:adm3}), when (\ref{i:adm1}) and (\ref{i:adm2}) hold. 
We only have to show that each $\mathbf{u}$ satisfying the conditions of~(\ref{i:adm3}) starts with~$01$ (and, symmetrically, $\mathbf{v}$~starts with~$10$). 
If $\overline{\mathbf{u}}$ started with~$00$, then we had $\overline{\mathbf{u}} = \overline{0}$ and thus $g\big(\Omega^{(\overline{\mathbf{u}},\overline{\mathbf{v}})}\big) = 1$, contradicting the assumption that $g\big(\Omega^{(\overline{\mathbf{u}},\overline{\mathbf{v}})}\big) = g\big(\Omega^{(\mathbf{a},\mathbf{b})}\big) > 1$. 

The converse implication is a direct consequence of Lemma~\ref{l:entpos}.

\subsection{Alternating lexicographic order}
For the alternating case, the following lemma and the subsequent remarks show that the condition $g\big(\Omega^{(\overline{\mathbf{u}},\overline{\mathbf{v}})}\big) = g\big(\Omega^{(\mathbf{a},\mathbf{b})}\big)$ is often easy to verify.
(Recall that $g\big(\Omega^{(\overline{\mathbf{u}},\overline{\mathbf{v}})}\big)$ can be determined by Lemma~\ref{l:K}.)
However, Example~\ref{ex:3} in Section~\ref{sec:main-results} shows that, contrary to the case of positive~$\beta$, it is not sufficient that $\mathbf{a}, \mathbf{b} \in \{\mathbf{u}, \mathbf{v}\}^\omega$ and $(\overline{\mathbf{u}}; \overline{\mathbf{v}})$ is a pair of critical itineraries.

\begin{lemma} \label{l:entneg}
Let $\mathbf{a} \in 0\, \Omega$, $\mathbf{b} \in 1\, \Omega$, $\mathbf{u} \in 0\, \{0,1\}^*$, $\mathbf{v} \in 1\, \{0,1\}^*$, such that $\mathbf{a}, \mathbf{b} \in \{\mathbf{u}, \mathbf{v}\}^\omega$, $\mathbf{a} \in \Omega^{(\mathbf{a},\mathbf{b}]}$, $\mathbf{b} \in \Omega^{[\mathbf{a},\mathbf{b})}$, $\overline{\mathbf{u}} \in \Omega^{(\overline{\mathbf{u}},\overline{\mathbf{v}}]}$, and $\overline{\mathbf{v}} \in \Omega^{[\overline{\mathbf{u}},\overline{\mathbf{v}})}$, where $\Omega$ is equipped with the alternating lexicographic order.
If $g(\Omega^{(\overline{\mathbf{u}},\overline{\mathbf{v}})})^{-|\mathbf{u}|} + g(\Omega^{(\overline{\mathbf{u}},\overline{\mathbf{v}})})^{-|\mathbf{v}|} < 1$, then
\[
g\big(\Omega^{(\overline{\mathbf{u}},\overline{\mathbf{v}})}\big)  = g\big(\Omega^{(\mathbf{a},\mathbf{b})}\big).
\]
\end{lemma}

\begin{proof}
Let 
\begin{align*}
\mathbf{a}_{\min} & = \overline{\mathbf{u}}, & \hspace{-.5em} \mathbf{a}_{\max} & = \mathbf{u} \overline{\mathbf{v}}, & \hspace{-.5em} \mathbf{b}_{\min} & = \mathbf{v} \overline{\mathbf{u}}, & \hspace{-.5em} \mathbf{b}_{\max} & = \overline{\mathbf{v}}, & & \hspace{-.5em} \text{if $|\mathbf{u}|$ and $|\mathbf{v}|$ are even}, \\
\mathbf{a}_{\min} & = \overline{\mathbf{uv}}, & \hspace{-.5em} \mathbf{a}_{\max} & = \mathbf{u} \overline{\mathbf{uv}}, & \hspace{-.5em} \mathbf{b}_{\min} & = \mathbf{v} \overline{\mathbf{vu}}, & \hspace{-.5em} \mathbf{b}_{\max} & = \overline{\mathbf{vu}}, & & \hspace{-.5em} \text{if $|\mathbf{u}|$ and $|\mathbf{v}|$ are odd}, \\
\mathbf{a}_{\min} & = \overline{\mathbf{u}}, & \hspace{-.5em} \mathbf{a}_{\max} & = \mathbf{uv} \overline{\mathbf{u}}, & \hspace{-.5em} \mathbf{b}_{\min} & = \mathbf{vv} \overline{\mathbf{u}}, & \hspace{-.5em} \mathbf{b}_{\max} & = \mathbf{v} \overline{\mathbf{u}}, & & \hspace{-.5em} \text{if $|\mathbf{u}|$ is even and $|\mathbf{v}|$ is odd}, \\
\mathbf{a}_{\min} & = \mathbf{u} \overline{\mathbf{v}}, & \hspace{-.5em} \mathbf{a}_{\max} & = \mathbf{uu} \overline{\mathbf{v}}, & \hspace{-.5em} \mathbf{b}_{\min} & = \mathbf{vu} \overline{\mathbf{v}}, & \hspace{-.5em} \mathbf{b}_{\max} & = \overline{\mathbf{v}}, & & \hspace{-.5em} \text{if $|\mathbf{u}|$ is odd and $|\mathbf{v}|$ is even}.
\end{align*}
Then $\mathbf{a}_{\min} \le \mathbf{a} \le \mathbf{a}_{\max}$ and $\mathbf{b}_{\min} \le \mathbf{b} \le \mathbf{b}_{\max}$, thus
\[
\Omega^{(\mathbf{a}_{\min},\mathbf{b}_{\max})} \subseteq \Omega^{(\mathbf{a},\mathbf{b})} \subseteq \Omega^{(\mathbf{a}_{\max},\mathbf{b}_{\min})}
\]
and $\Omega^{(\mathbf{a}_{\min},\mathbf{b}_{\max})} \subseteq \Omega^{(\overline{\mathbf{u}},\overline{\mathbf{v}})} \subseteq \Omega^{(\mathbf{a}_{\max},\mathbf{b}_{\min})}$.
We first prove that
\begin{equation} \label{e:capuvneg}
\Omega^{(\mathbf{a}_{\max},\mathbf{b}_{\min})} \cap [\mathbf{a}_{\min},\mathbf{b}_{\max}] \subseteq \{\mathbf{u}, \mathbf{v}\}^\omega.
\end{equation}
Let $\mathbf{c} \in \Omega^{(\mathbf{a}_{\max},\mathbf{b}_{\min})} \cap [\mathbf{a}_{\min},\mathbf{b}_{\max}]$, and consider an arbitrary decomposition $\mathbf{c} = \mathbf{wd}$ with $\mathbf{w} \in \{\mathbf{u}, \mathbf{v}\}^*$, $\mathbf{d} \in \Omega$.
Assume that $\mathbf{d}$ starts with~$0$.
Then $\mathbf{d}$ starts with $\mathbf{u}$ because $\mathbf{a}_{\min} \le \mathbf{d} \le \mathbf{a}$.
Here, the inequality $\mathbf{a}_{\min} \le \mathbf{d}$ comes from the following considerations.
If $|\mathbf{u}|$ and~$|\mathbf{v}|$ are even, then $\mathbf{d} \ge \overline{\mathbf{u}}$ as in the proof of Lemma~\ref{l:entpos}.
If $|\mathbf{u}|$ and~$|\mathbf{v}|$ are odd, then $\mathbf{d} \ge \overline{\mathbf{uv}}$ follows from $\mathbf{c} \in [\overline{\mathbf{uv}}, \overline{\mathbf{vu}}]$ when $\mathbf{w} \in \{\mathbf{uv}\}^* \cup \mathbf{v}\, \{\mathbf{uv}\}^*$ and from $\mathbf{c} \in \Omega^{(\mathbf{u}\overline{\mathbf{uv}},\mathbf{v}\overline{\mathbf{vu}})}$ when $\mathbf{w}$ ends with a word in $\{\mathbf{u}, \mathbf{vv}\}\, \{\mathbf{uv}\}^*$.
If $|\mathbf{u}|$ is even and $|\mathbf{v}|$ is odd, then $\mathbf{d} \ge \overline{\mathbf{u}}$ follows from $\mathbf{c} \in [\overline{\mathbf{u}}, \mathbf{v} \overline{\mathbf{u}}]$ when $\mathbf{w} \in \{\mathbf{u}\}^* \cup \mathbf{v}\, \{\mathbf{u}\}^*$ and from $\mathbf{c} \in \Omega^{(\mathbf{uv}\overline{\mathbf{u}},\mathbf{vv}\overline{\mathbf{u}})}$ when $\mathbf{w}$ ends with a word in $\{\mathbf{uv}, \mathbf{vv}\}\, \{\mathbf{u}\}^*$.
Finally, if $|\mathbf{u}|$ is odd and $|\mathbf{v}|$ is even, then $\mathbf{d} \ge \mathbf{u} \overline{\mathbf{v}}$ follows from $\mathbf{c} \in [\mathbf{u} \overline{\mathbf{v}}, \overline{\mathbf{v}}]$ when $\mathbf{w}$ is the empty word and from $\mathbf{c} \in \Omega^{(\mathbf{uu}\overline{\mathbf{v}},\mathbf{vu}\overline{\mathbf{v}})}$ when $\mathbf{w}$ ends with $\mathbf{u}$ or~$\mathbf{v}$.
Symmetrically, we obtain that  $\mathbf{d}$ starts with~$\mathbf{v}$ whenever it starts with~$1$. 
Therefore, \eqref{e:capuvneg} holds, and we have
\[
\Omega^{(\mathbf{a}_{\min},\mathbf{b}_{\max})} \subseteq \Omega^{(\mathbf{a}_{\max},\mathbf{b}_{\min})} \subseteq \Omega^{(\mathbf{a}_{\min},\mathbf{b}_{\max})} \cup \bigcup_{n=0}^\infty \Omega^{(\mathbf{a}_{\min},\mathbf{b}_{\max})}_n\, \{\mathbf{u}, \mathbf{v}\}^\omega,
\]
thus $g(\Omega^{(\mathbf{a}_{\max},\mathbf{b}_{\min})}) = \max(g(\Omega^{(\mathbf{a}_{\min},\mathbf{b}_{\max}})), g\big(\bigcup_{n=0}^\infty S^n \{\mathbf{u}, \mathbf{v}\}^\omega)\big)$.

Since $g\big(\bigcup_{n=0}^\infty S^n \{\mathbf{u}, \mathbf{v}\}^\omega\big)\big)$ is the only solution of $x^{-|\mathbf{u}|} + x^{-|\mathbf{v}|} = 1$ with $x \ge 1$, $g(\Omega^{(\overline{\mathbf{u}},\overline{\mathbf{v}})})^{-|\mathbf{u}|} + g(\Omega^{(\overline{\mathbf{u}},\overline{\mathbf{v}})})^{-|\mathbf{v}|} < 1$ implies that $g\big(\bigcup_{n=0}^\infty S^n \{\mathbf{u}, \mathbf{v}\}^\omega\big)\big) < g\big(\Omega^{(\overline{\mathbf{u}},\overline{\mathbf{v}})}\big) \le g\big(\Omega^{(\mathbf{a}_{\max},\mathbf{b}_{\min})}\big)$, thus  $g\big(\Omega^{(\mathbf{a}_{\max},\mathbf{b}_{\min})}\big) = g\big(\Omega^{(\mathbf{a}_{\min},\mathbf{b}_{\max}})\big)$, which gives that $g\big(\Omega^{(\overline{\mathbf{u}},\overline{\mathbf{v}})}\big) = g\big(\Omega^{(\mathbf{a},\mathbf{b})}\big)$.
\end{proof}

From the proof of Lemma~\ref{l:entneg}, we can also derive other conditions that guarantee $g\big(\Omega^{(\overline{\mathbf{u}},\overline{\mathbf{v}})}\big) = g\big(\Omega^{(\mathbf{a},\mathbf{b})}\big)$, e.g., $g(\Omega^{(\mathbf{a}_{\min},\mathbf{b}_{\max}}))^{-|\mathbf{u}|} + g(\Omega^{(\mathbf{a}_{\min},\mathbf{b}_{\max}}))^{-|\mathbf{v}|} \le 1$, or $g(\Omega^{(\overline{\mathbf{u}},\overline{\mathbf{v}})})^{-|\mathbf{u}|} + g(\Omega^{(\overline{\mathbf{u}},\overline{\mathbf{v}})})^{-|\mathbf{v}|} \le 1$ and $g\big(\Omega^{(\overline{\mathbf{u}},\overline{\mathbf{v}})}\big) \le g\big(\Omega^{(\mathbf{a},\mathbf{b})}\big)$.

\section{Lorenz Maps}  \label{sec:G} 
A~\textit{Lorenz map}, as defined, e.g., in~\cite{HS}, is a function $f:\, [0,1] \to [0,1]$ satisfying:
\begin{enumerate}
\item 
There exists a $c \in (0,1)$ such that $f$ is continuous and strictly increasing on $[0,c)$ and on $(c,1]$;
\item 
$\lim_{x\uparrow c} f(x) =1$ and $\lim_{x\downarrow c} f(x) =0$.
\end{enumerate}
For $\beta > 1$, $1 \le p \le \frac{1}{\beta-1}$, the restriction of $f_{\beta,p}$ to $\big[0, \frac{\beta}{\beta-1}\big]$ is thus conjugate to the Lorenz map with constant slope~$\beta$ and $c = \frac{\beta-1}{\beta}\, p$. 

In~\cite{HS}, the authors define a fairly weak notion of what it means for a Lorenz map to be expanding. 
Specifically, a Lorenz map is said to be \textit{topologically expansive} if there exists an $\epsilon>0$ such that any two distinct forward orbits $(x_{0},x_{1},x_{2}, \dots)$ and $(y_{0},y_{1},y_{2}, \dots)$ satisfy $|x_{i}-y_{i}|\geq\epsilon$ for some $i \ge 0$. 
They prove that a pair $(\mathbf{a}; \mathbf{b})$ of binary strings is a pair of critical itineraries of a topologically expansive Lorenz map if and only if $(\mathbf{a}; \mathbf{b})$ satisfy condition~(\ref{i:adm1}) in Definition~\ref{def:admissible}.

In~\cite{Gl,GS}, for example, the authors define a stronger notion of what it means for a Lorenz map to be expanding. 
Specifically, a Lorenz map is an $L_{1+\epsilon}$ map if $f$ is differentiable except at the point $c$ of discontinuity and if there exists an $\epsilon>0$ such that $f'(x) \geq 1+\epsilon$ for all $x\neq c$.  
In \cite{P2}, the author proves that if an $L_{1+\epsilon}$ map $f$ is transitive, then $f$ is topologically conjugate to a (generalized) $\beta$-transformation.  
In~\cite{GS}, the authors give necessary and sufficient conditions for an $L_{1+\epsilon}$ map to be transitive in terms of its critical itineraries. 
In~\cite{Gl}, the author gives necessary and sufficient conditions for an $L_{1+\epsilon}$ map to be topologically conjugate to a generalized $\beta$-transformation in terms of its critical itineraries.
There does not, however, seem to be a characterization of the pairs of
critical itineraries of an $L_{1+\epsilon}$ map.
The pair $(\mathbf{a}; \mathbf{b})$ with
\[
\mathbf{a} =011\overline{100}, \qquad \mathbf{b}=100\overline{011},
\]
is the pair of critical itineraries of an $L_{1+\epsilon}$ map, but is not, according to the main result of this paper, the pair of critical itineraries of a generalized $\beta$-transformation.
The pair $(\mathbf{a}; \mathbf{b})$ with 
\[
\mathbf{a}=  0111\overline{10}, \qquad \mathbf{b}=\overline{10011110},
\]
is an example of a pair of critical itineraries of a topologically expanding Lorenz map, but not the pair of critical itineraries of an $L_{1+\epsilon}$ Lorenz map. 
This can be verified by noting that, if  $(\mathbf{a}; \mathbf{b})$  were the pair of critical itineraries of an $L_{1+\epsilon}$ map, then by the criteria in~\cite{Gl}, it would be the pair of critical itineraries of a generalized $\beta$-transformation. 
But, by the main theorem of this paper, that is not the case.

\bibliographystyle{amsplain}
\bibliography{critical}
\end{document}